\documentclass{article}

\usepackage{arxiv}

\usepackage[utf8]{inputenc} 
\usepackage[T1]{fontenc}    
\usepackage{hyperref}       
\usepackage{url}            
\usepackage{booktabs}       
\usepackage{amsfonts}       
\usepackage{nicefrac}       
\usepackage{microtype}      
\usepackage{lipsum}		
\usepackage{graphicx}
\usepackage[authoryear]{natbib}
\usepackage{doi}

\usepackage{mathtools,upgreek,eucal,mathrsfs,enumitem,amsthm,amsmath}
\usepackage{times,epsfig,makeidx,amsfonts,amsmath,dcolumn,multirow,amssymb,colortbl,bm,url}
\usepackage{algorithm}
\usepackage[]{algorithmic}
 \usepackage[flushleft]{threeparttable}
\usepackage{subfig}

\newtheorem{theorem}{Theorem}
\newtheorem{definition}{Definition}

\newcommand{\bxi}{\bm{\bxi}}

\newcommand{\btheta}{\bm{\theta}}

\newcommand{\blambda}{\bm{\lambda}}

\newcommand{\bbeta}{\bm{\beta}}

\newcommand{\bvarphi}{\bm{\varphi}}

\newcommand{\bJ}{{\bf J}}

\newcommand{\bX}{{\bf X}}
\newcommand{\bx}{{\bf x}}
\newcommand{\by}{{\bf y}}

\newcommand{\bone}{{\bf 1}}

\newcommand{\bbexp}{{\mathbb E}}
\newcommand{\R}{{\mathbb R}}

\title{On a prior based on the Wasserstein information matrix}

\author{ \href{https://orcid.org/0000-0002-2218-5734}{\includegraphics[scale=0.06]{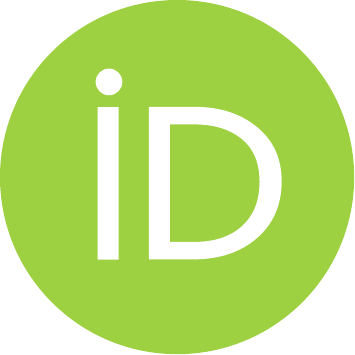}\hspace{1mm}Wuchen Li}\\
Department of Mathematics\\
University of South Carolina\\
South Carolina, USA.\\
	\texttt{wuchen@mailbox.sc.edu} \\
	\And
	\href{https://orcid.org/0000-0001-7183-8407}{\includegraphics[scale=0.06]{orcid.pdf}\hspace{1mm}Francisco Javier Rubio} \\
	Department of Statistical Science\\
	University College London \\
	London, UK\\
	\texttt{f.j.rubio@ucl.ac.uk} 
	}

\hypersetup{
pdftitle={XXX},
pdfsubject={stat.ME, stat.AP},
pdfauthor={Li, Rubio},
}

\begin{document}
\maketitle

\begin{abstract}
We introduce a prior for the parameters of univariate continuous distributions, based on the Wasserstein information matrix, which is invariant under reparameterisations. 
We discuss the links between the proposed prior with information geometry. 
We present sufficient conditions for the propriety of the posterior distribution for general classes of models. We present a simulation study that shows that the induced posteriors have good frequentist properties.
\end{abstract}



\section{Introduction}\label{sec:intro}
In Bayesian parametric inference, the choice of the prior plays a fundamental role. In scenarios where the prior information about the model parameters is vague or unreliable, it is desirable to use priors which do not require the user to specify their parameters (hyperparameters). The main aim of \emph{Objective Bayes} \citep{berger:2006,consonni:2018} is indeed to produce priors via formal rules \citep{kass:1996}, which typically depend only on the statistical model. Such rules usually aim at producing a prior that has little effect on the inference on the parameters, or that is invariant under reparameterisations, or that penalises the model complexity. Priors obtained with formal rules are usually referred to as \emph{Objective priors} or \emph{Non-informative priors}. 
We refer the reader to \cite{leisen:2020} for a recent review of methods for constructing priors based on formal rules.
A pioneering contribution in this area is the \emph{Jeffreys prior} \citep{jeffreys:1946}, which is obtained by calculating the square root of the determinant of the Fisher information matrix (FIM) \citep{robert:2009}. The aim behind the construction of the Jeffreys prior is to produce a prior that is invariant under reparameterisations.

Another direction for constructing a prior based on a formal rule consists of looking at the genesis of the Jeffreys prior \citep{kass:1996}. The Jeffreys prior is typically motivated by its invariance under reparameterisations, however, it can also be motivated using concepts from information geometry \citep{amari:2016,nielsen:2020,amari:2021b,amari:2022}. Briefly, the Kullback-Leibler divergence behaves locally as a function of a distance function determined by the Riemannian metric. The Jeffreys prior can be seen as the natural volume associated to such metric, and natural volume elements generate uniform measures on manifolds \citep{kass:1996}. Moreover, natural volumes of Riemannian metrics are invariant under reparameterisations \citep{kass:1996}.
Intuitively, this suggests that other distances could be used to construct alternative priors. 
In this line, a natural alternative consists of using the optimal transport induced information matrix \citep{li:2019}, referred to as the Wasserstein information matrix (WIM). 
The construction of the WIM can be justified using ideas from ``transport information geometry'', which is the intersection between optimal transport \citep{villani:2003} and information geometry \citep{amari:2016,amari:2021b}. We refer the reader to \cite{Li2021_transportd,Li2021_transportc} and \cite{amari:2021b} for a more extensive treatment of this area. 
The idea behind the construction of the WIM consists of using tools from optimal transport, where a distance between distributions is used to construct an information matrix. 
\cite{li:2019} focused on the particular choice of the Wasserstein-2 distance. The Wasserstein-2 distance can be associated to a metric operator, namely the WIM, which is different in nature from the Fisher information matrix \citep{amari:2021b}. 
Although a vast amount of literature has been devoted to the study of the Fisher information matrix and the Jeffreys prior, there is a void in the study of priors associated to the Wasserstein information matrix.  

We propose a formal rule for constructing a prior, which is invariant under reparameterisations, based on the Wasserstein information matrix. 
The construction of this prior (referred to as the Wasserstein prior hereafter) is analogous to that of the Jeffreys prior. However, as shown later, we find that the Wasserstein prior has a different functional form for several models and, appealingly, requires a lower order of differentiability. This helps overcome some challenges with the Jeffreys prior, where the required higher order of differentiability precludes its construction for some non-regular models \citep{shemyakin:2014,li:2019}. Moreover, as we will show later in the simulation study, the Wasserstein prior induces a posterior with good frequentist properties in the models studied here. 


\section{Wasserstein information matrix}\label{sec:wim}
Let $X$ be a continuous random variable with finite second moment, and $F(x\mid\btheta)$ be the corresponding cumulative distribution function (cdf) with support $\mathcal{D} \subset {\mathbb R}$, and parameter $\btheta \in \Theta \subset {\mathbb R}^d$, with $d\geq 1$. Let us assume that $F(x\mid\btheta)$ is absolutely continuous, and let $f(x\mid\btheta)$ be the corresponding probability density function (pdf). 

Consider the Wasserstein information matrix (WIM) proposed in \cite{li:2019}
\begin{eqnarray}\label{eq:WIM}
W_{ij}(\btheta) &=& \bbexp \left[ \dfrac{\dfrac{\partial}{\partial \btheta_i}F(X\mid\btheta)  \dfrac{\partial}{\partial \btheta_j}F(X\mid\btheta)}{f(X\mid\btheta)^2} \right],
\end{eqnarray}
where the expectation is taken with respect to $F(x\mid\btheta)$. A clear difference between the WIM and the FIM is that the former is based on derivatives of the cdf (with respect to the parameters), while the latter is based on derivatives of the pdf. This is an appealing property as it reduces the conditions for the existence of the WIM \citep{li:2019}, allowing its construction for non-regular models.
Next, we present a brief description of the motivation behind the construction of the WIM. The details are somewhat technical, but we refer the reader to  \cite{li:2019} for a detailed derivation of the WIM.

As discussed in Section \ref{sec:intro}, a distance between probability distributions can be used to define an information matrix. In our case, we focus on the analysis of the information matrix (WIM) implied by the Wasserstein-2 distance. Given two parameter values $\btheta_0$, $\btheta_1\in\Theta$, the Wasserstein-2 distance between two probability distributions with support on $\mathcal{D} \subset {\mathbb R}$, $F(\cdot \mid\btheta_0)$ and $F(\cdot \mid\btheta_1)$, satisfies the following relationship with the corresponding quantile functions \citep{villani:2003}
\begin{equation*}
\mathrm{Dist}_{W}(F(\cdot \mid \btheta_0), F(\cdot \mid \btheta_1)) = \sqrt{ \int_{0}^{1} \left| F^{-1}(u\mid \btheta_0) - F^{-1}(u\mid\btheta_1) \right|^2 du},
\end{equation*}
where $F^{-1}$ is the quantile function associated to the cdf $F$. It can be shown that the Wasserstein-2 distance $\mathrm{Dist}_W$ defines a Riemannian metric among probability distributions \citep{villani:2003}, which can be used to establish the connection between such metric with an information matrix. More specifically, the infinitesimal expansion of the squared Wasserstein-2 distance establishes a link of this metric and the Wasserstein information matrix. That is, let $\Delta\btheta=(\Delta\btheta_1,\dots,\Delta \btheta_d)^{\top}\in\mathbb{R}^d$ such that $\btheta_0+\Delta\btheta \in \Theta$, one can show that \citep{li:2019}
\begin{equation*}
\mathrm{Dist}_W(F(\cdot \mid \btheta_0), F(\cdot \mid \btheta_0+\Delta\btheta))^2=\sum_{i,j=1}^d W_{ij}(\btheta_0)\Delta\btheta_i\Delta\btheta_j+o(\|\Delta\btheta\|^2).     
\end{equation*}
This shows a link between the Wasserstein-2 distance and the WIM, which is discussed in detail in section 7.7 of \cite{amari:2021b}. We can also see from this result that the WIM shares a similar derivation to that of the Fisher information matrix (see Chapter 5 of \citep{ghosh:2006} for an extensive discussion). We remark that one can also define the WIM in higher dimensional sample spaces (that is, for random vectors). However, this requires solving an elliptic partial differential equation \citep{li:2019}.


\section{The Wasserstein prior}\label{sec:wprior}
In this section, we propose the Wasserstein prior, whose main motivation is to obtain an invariant prior. We also describe the precise meaning of the invariance property and its connection with the Jeffreys prior.
\subsection{One parameter case}
Consider the case where $d=1$, that is, we focus on the case where $F(x\mid \theta)$ contains only one parameter. Then, the WIM \eqref{eq:WIM} becomes  
\begin{eqnarray*}\label{eq:WIM1D}
W(\theta) &=& \bbexp\left[ \dfrac{\left\{\dfrac{\partial}{\partial \theta}F(X\mid\theta)\right\}^2 }{f(X\mid\theta)^2} \right].
\end{eqnarray*}
Let $\varphi = h(\theta)$ be a reparameterisation of $F(x\mid\theta)$. Let us denote the WIM associated to $F(x\mid\theta)$ by $W(\theta)$, and the WIM associated to $F(x\mid\varphi)$ by $\tilde{W}(\varphi)$. From the above expression, we can see that 
\begin{eqnarray*}
\tilde{W}(\varphi) = W(\theta)\left(\dfrac{d \theta}{d \varphi}\right)^2 .
\end{eqnarray*}
Indeed, the FIM also satisfies this relationship \citep{robert:2009}. This suggests the construction of an invariant prior, based on the WIM, in a similar fashion as the Jeffreys prior (which is based on the FIM). Define the prior (up to a positive proportionality constant)
\begin{eqnarray*}
\pi_W(\theta) \propto \sqrt{W(\theta)}.
\end{eqnarray*}
It follows that this prior is invariant under reparameterisations in the sense that
\begin{eqnarray*}
{\pi}_{\tilde{W}}(\varphi)= \pi_W(\theta) \left \vert \dfrac{d \theta}{d\varphi}  \right \vert,
\end{eqnarray*}
where ${\pi}_{\tilde{W}}(\varphi) \propto \sqrt{\tilde{W}(\varphi)}$. That is, the priors $\pi_W(\theta)$ and ${\pi}_{\tilde{W}}(\varphi)$ are related by the corresponding change of variable.
Therefore, this represents a strategy for constructing a prior based on a formal rule \citep{kass:1996} which is invariant under reparameterisations, in the same spirit as the invariance property of the Jeffreys prior \citep{jeffreys:1946}. We formalise this idea next.

\subsection{Multi-parameter case}

Consider now the general case $d\geq1$ and let $\bvarphi = h(\btheta)$ be a reparameterisation of $F(x\mid\btheta)$. Note first that the WIM of $\bvarphi$, $\tilde{W}(\bvarphi)$, can be written after a change of variable as:
\begin{eqnarray*}
\tilde{W}(\bvarphi) = \bJ^{\top} W(\btheta) \bJ,
\end{eqnarray*}
where $\bJ$ is the Jacobian matrix with entries
\begin{eqnarray*}
\bJ_{ij} = \dfrac{\partial \btheta_i}{\partial \bvarphi_j}.
\end{eqnarray*}
The proof of this result is analogous to the proof of the invariance property of the FIM, which can be found in \cite{lehmann:2006}.
Consequently, we have that
\begin{eqnarray*}
\mbox{det} \, \tilde{W}(\bvarphi) = \mbox{det} \, W(\btheta) (\mbox{det} \, \bJ )^2.
\end{eqnarray*}
This result suggests the construction of an invariant prior, based on the WIM, in a similar fashion as the Jeffreys prior is obtained from the FIM. The construction of this prior is formalised in the following definition. 
\begin{definition}
The Wasserstein prior is defined, up to a positive proportionality constant, as
\begin{eqnarray}\label{eq:piw}
\pi_W(\btheta) \propto \sqrt{\mbox{det} \, W(\btheta)},
\end{eqnarray}
where $W(\btheta)$ denotes the Wasserstein information matrix \eqref{eq:WIM}.
\end{definition}


\section{Examples}\label{sec:examples}
In this section, we present three examples where we illustrate the calculation of the WIM and the Wasserstein prior. In all cases, we provide sufficient conditions for the propriety of the posterior distribution.

\subsection*{The location-scale family}
Let $f_0$ be a symmetric and unimodal pdf with mode at $0$ and support on ${\mathbb R}$, and $F_0$ be the corresponding cdf. Let 
\begin{eqnarray}\label{eq:locsc}
F(x \mid \mu,\sigma) &=& F_0\left(\dfrac{x-\mu}{\sigma}\right),\quad f(x \mid \mu,\sigma) = \dfrac{1}{\sigma}f\left(\dfrac{x-\mu}{\sigma}\right), \quad x\in{\mathbb R},
\end{eqnarray}
denote the cdf and pdf of the class of symmetric and unimodal location-scale family of distributions, with location parameter $\mu\in{\mathbb R}$ and scale parameter $\sigma\in {\mathbb R}_+$. 

\begin{theorem}\label{th:locsc}
Suppose that $\int_{-\infty}^{\infty} t^2f_0(t)dt < \infty$.
The WIM and the Wasserstein prior of $(\mu,\sigma)$ in the location-scale family \eqref{eq:locsc} are:
\begin{eqnarray}\label{eq:piwls}
W(\mu,\sigma) = \left(
\begin{array}{c c}
1 & 0 \\
0 &  \int_{-\infty}^{\infty} t^2 f_0(t)dt.
\end{array}
\right), \quad
\pi_W(\mu,\sigma) \propto 1.
\end{eqnarray}
\end{theorem}

An important class of location-scale models is the family of scale mixtures of normal distributions. A pdf $f_0$ is said to belong to family of scale mixtures of normal distributions if it can be represented as:
\begin{eqnarray}\label{eq:SMN}
f_0(x) = \int_0^\infty \lambda^{\frac{1}{2}} \exp\left\{ - \dfrac{\lambda x^2}{2} \right\} d H(\lambda),
\end{eqnarray}
where $H$ is a cumulative distribution function with support on ${\mathbb R}_+$. The family of scale mixtures of normal distributions contains important distributions such as the Normal, Logistic, Laplace, Student-$t$, among other distributions (see \citealp{rubio:2014} for a discussion). The next result provides sufficient conditions for the propriety of the posterior distribution of $(\mu,\sigma)$ under the Wasserstein prior \eqref{eq:piwls} for the case when $f_0$ belongs to the family of scale mixtures of normal distributions. 

\begin{theorem}\label{th:properls}
Let $\bx = (x_1,\dots, x_n)^{\top}$ be an \textit{i.i.d.}~sample from \eqref{eq:locsc} with $f_0$ given by \eqref{eq:SMN}. Suppose that $\int_{-\infty}^{\infty} t^2f_0(t)dt < \infty$. Then, the posterior distribution of $(\mu,\sigma)$ associated to the Wasserstein prior \eqref{eq:piwls} is proper if $n > 2$ and 
$$\int_0^{\infty} \lambda^{1/2} dH(\lambda) < \infty.$$
\end{theorem}
\subsection*{The skew-normal distribution}
We now present a result in a one-parameter model, where we obtain the Wasserstein prior for the skewness parameter of the skew-normal distribution \citep{azzalini:1985}. Let $ \phi(x)$ and $ \Phi(x)$ be the pdf and cdf of the standard normal distribution. The skew-normal pdf is defined as \citep{azzalini:1985}:
\begin{eqnarray}\label{eq:SN}
s(x \mid \alpha) = 2 \phi(x) \Phi(\alpha x), \quad x\in{\mathbb R},
\end{eqnarray}
where $\alpha \in {\mathbb R}$ is a skewness parameter. The following result characterises the WIM and Wasserstein prior of $\alpha$.

\begin{figure}[h!]
\begin{tabular}{c c}
\includegraphics[scale=0.4]{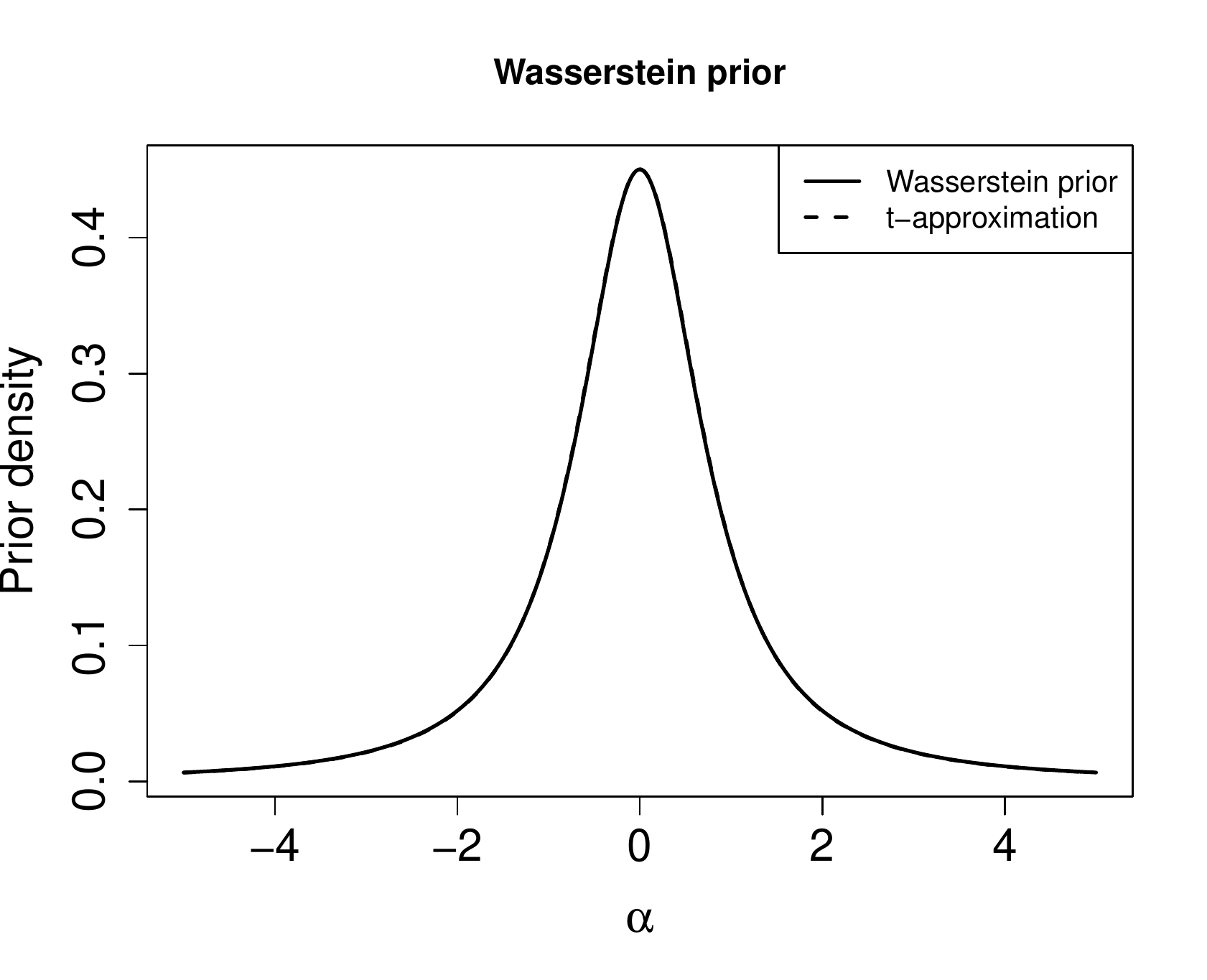} & \includegraphics[scale=0.4]{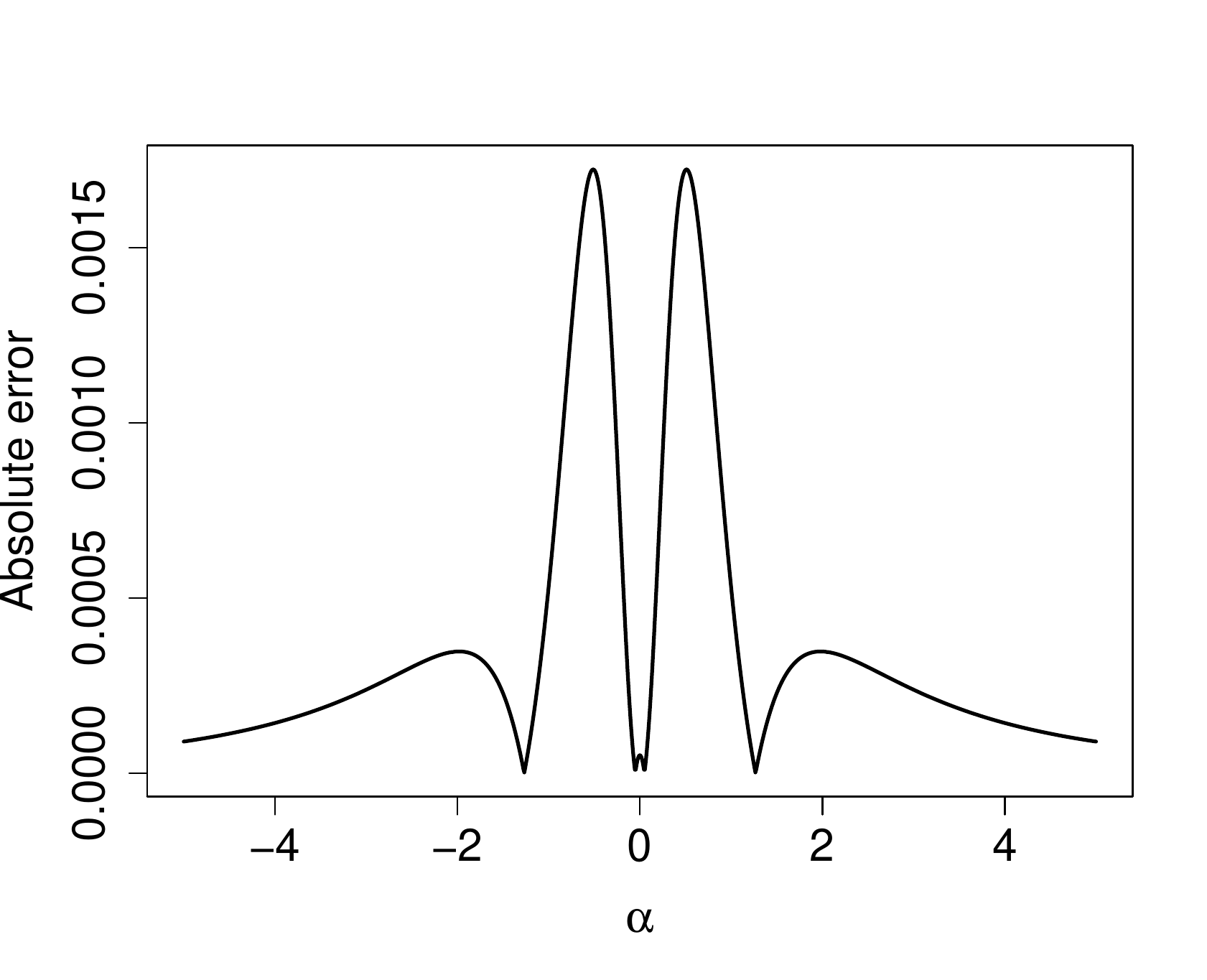}\\
(a) & (b)
\end{tabular}
\caption{(a) Wasserstein prior and Student-$t$ approximation. (b) Absolute error of the Student-$t$ approximation.}
\label{fig:tapp}
\end{figure}

\begin{theorem}\label{th:sn}
Consider the skew-normal distribution \eqref{eq:SN}. Then,
\begin{itemize}
\item[(i)] The WIM of $\alpha$ is given by
\begin{eqnarray*}
W(\alpha) = \int_{-\infty}^{\infty} \frac{\sqrt{2} e^{-\frac{1}{2} \left(2 \alpha ^2+1\right) x^2}}{\pi ^{3/2} \left(\alpha
   ^2+1\right)^2 \left(\mbox{erf}\left(\frac{\alpha  x}{\sqrt{2}}\right)+1\right)}dx.
\end{eqnarray*}
\item[(ii)] The Wasserstein prior
\begin{eqnarray}\label{eq:piwsn}
\pi_W(\alpha)  \propto \sqrt{\int_{-\infty}^{\infty} \frac{\sqrt{2} e^{-\frac{1}{2} \left(2 \alpha ^2+1\right) x^2}}{\pi ^{3/2} \left(\alpha
   ^2+1\right)^2 \left(\mbox{erf}\left(\frac{\alpha  x}{\sqrt{2}}\right)+1\right)}dx},
\end{eqnarray}
is symmetric about $0$.
\item[(iii)]   $\pi_W(\alpha)$ is integrable.
\item[(iv)] The tails of $\pi_W(\alpha)$ are of order ${\mathcal O}(\vert \alpha \vert^{-5/2})$.
\end{itemize}
\end{theorem}
The tail behaviour of the Wasserstein prior $\pi_W(\alpha)$ differs from that of the Jeffreys prior of $\alpha$ \citep{rubio:2014b}, which has tails of order ${\mathcal O}(\vert \alpha \vert^{-3/2})$, and the total variation prior proposed in \cite{dette:2018}, which has tails of order ${\mathcal O}(\vert \alpha \vert^{-2})$.
The characterisation of the propriety and tail behaviour of $\pi_W(\alpha)$ in the previous theorem suggests that one could approximate it using a symmetric distribution with the same tail behaviour. A natural candidate is the Student-$t$ distribution with $\nu = 3/2$ degrees of freedom. We found that a scale parameter $\sigma_t = 0.757$ produces a good approximation in the main body of the distribution, while the tails have the exact same weight (see Figure \ref{fig:tapp}). 

In the next theorem, we construct a prior for the skew normal distribution with location and scale parameters $(\mu,\sigma)$ and skewness parameter $\alpha$, based on a product prior structure using the priors \eqref{eq:piwls} and \eqref{eq:piwsn}. We show that the posterior distribution is proper under mild conditions. This prior can be interpreted as an  \emph{Independence Wasserstein prior} (analogous to the \emph{independence Jeffreys prior}, \citealp{rubio:2014,rubio:2014b}), in the sense that it is constructed as the product of the Wasserstein priors for each parameter (or groups of parameters) while considering the other parameters as fixed.
\begin{theorem}\label{th:propersn}
Let $\bx=(x_1,\dots,x_n)^{\top}$ be an \textit{i.i.d.}~sample from the skew normal distribution with pdf
\begin{eqnarray*}
s(x \mid \mu,\sigma, \alpha) = \dfrac{2}{\sigma} \phi\left(\dfrac{x-\mu}{\sigma}\right) \Phi\left(\alpha \dfrac{x-\mu}{\sigma}\right).
\end{eqnarray*}
Consider the improper product prior structure using the Wasserstein priors \eqref{eq:piwls} and \eqref{eq:piwsn}
\begin{eqnarray}
\pi(\mu,\sigma,\alpha) \propto \sqrt{\int_{-\infty}^{\infty} \frac{\sqrt{2} e^{-\frac{1}{2} \left(2 \alpha ^2+1\right) x^2}}{\pi ^{3/2} \left(\alpha
   ^2+1\right)^2 \left(\mbox{erf}\left(\frac{\alpha  x}{\sqrt{2}}\right)+1\right)}dx}.
   \label{eq:iwpsn}
\end{eqnarray}
Then, the posterior distribution of $(\mu,\sigma,\alpha)$ is proper if $n > 2$.
\end{theorem}

\subsection*{Normal linear regression}
We now study the WIM and the Wasserstein prior for the normal linear regression model,
\begin{eqnarray}\label{eq:lrm}
y_i = \bx_i^{\top} \bbeta + \epsilon_i, \,\,\,\,\, i = 1,\dots,n.
\end{eqnarray}
where $\bx_i^{\top}\in{\mathbb R}^p$ is a vector of covariates, $\bbeta\in{\mathbb R}^p$ is a vector of regression coefficients, $\epsilon_i \stackrel{i.i.d.}{\sim} N(0,\sigma^2)$ denote the errors. Let $\bX = (\bx_1,\dots,\bx_n)^{\top}$ denote the design matrix and $\by = (y_1,\dots,y_n)^{\top}$ the vector of response variables.
\begin{theorem}\label{th:lrm}
Consider the linear regression model \eqref{eq:lrm}, and suppose that $\bX$ has full column rank. Then, the WIM and the Wasserstein prior are given by,
\begin{eqnarray}\label{eq:piwlrm}
W(\bbeta,\sigma) = \left(
\begin{array}{c c}
\bX^{\top}\bX & 0 \\
0 &  1
\end{array}
\right), \quad
\pi_W(\bbeta,\sigma) \propto 1.
\end{eqnarray}
\end{theorem}
The next result presents sufficient conditions for the propriety of the posterior distribution of $(\bbeta,\sigma)$ under the Wasserstein prior \eqref{eq:piwlrm}.
\begin{theorem}\label{th:properlrm}
Consider the Normal linear regression model \eqref{eq:lrm} together with the Wasserstein prior \eqref{eq:piwlrm}.  Suppose that $\bX$ has full column rank and that $\by$ is not in the column space of $\bX$. Then, the posterior distribution of $(\bbeta,\sigma)$ is proper if $n > p + 1$.
\end{theorem}

\section{Simulation Studies}\label{sec:simulation}
In this section we present two simulation studies to assess the performance of the posterior distributions induced by the Wasserstein prior. 

In the first simulation scenario, we evaluate the performance of the independence Wasserstein prior \eqref{eq:iwpsn} and compare it against the independence Jeffreys prior \citep{rubio:2014b}. We simulate $N=250$ samples of size $n=50,250,500$ from a skew-normal distribution \eqref{eq:SN} with $\mu=10$, $\sigma=1$, and $\lambda = 1,3,5$. We emphasise that the value $\lambda=1$ represents a very challenging scenario as the skew-normal distribution is weakly identifiable for $\vert \lambda \vert < 1.25$; in the sense that the skew-normal pdf is virtually symmetric for values of $\lambda$ in this region \citep{rubio:2016}.  In the second simulation scenario, we evaluate the performance of the Wasserstein prior \eqref{eq:lrm} in linear regression models. We simulate $N=250$ samples of size $n=50,250,500$ from the linear regression model \eqref{eq:lrm}, with $\bbeta = (1,0,0.5,1)^{\top}$ and $\sigma = 0.5$, where the first entry of $\bbeta$ represents the intercept. The entries of the design matrix are simulated from a multivariate normal distribution with zero mean, unit variance, and pairwise correlations of $0.5$. The values of $\bbeta$ are chosen to reflect different levels of signal-to-noise ratio and the effect of a spurious variable. For each of these samples, we simulate a posterior sample of size $1000$ using the R package `Rtwalk', using a burn-in period of $5000$ iterations and a thinning period of $25$ iterations (this is, a total of $300,000$ posterior samples were obtained for each sample). In all scenarios, we also compare the results against those associated to the maximum likelihood estimators (MLE). We choose the following performance measures to evaluate the different estimation methods and priors: `mMean' denoting the average of the posterior means across the $N$ simulated samples; `mSD' denoting the average of the posterior standard deviations; `mRMSE' denoting the average of the root mean squared errors; `Coverage' denoting the coverage proportion of the $95\%$ credible intervals; `mMLE' denoting the average of the maximum likelihood estimators; and `RMSE-MLE' denoting the root mean squared error of the maximum likelihood estimators across the $N$ samples.

\begin{figure}
\begin{tabular}{ccc}
\centering
\includegraphics[scale=0.5]{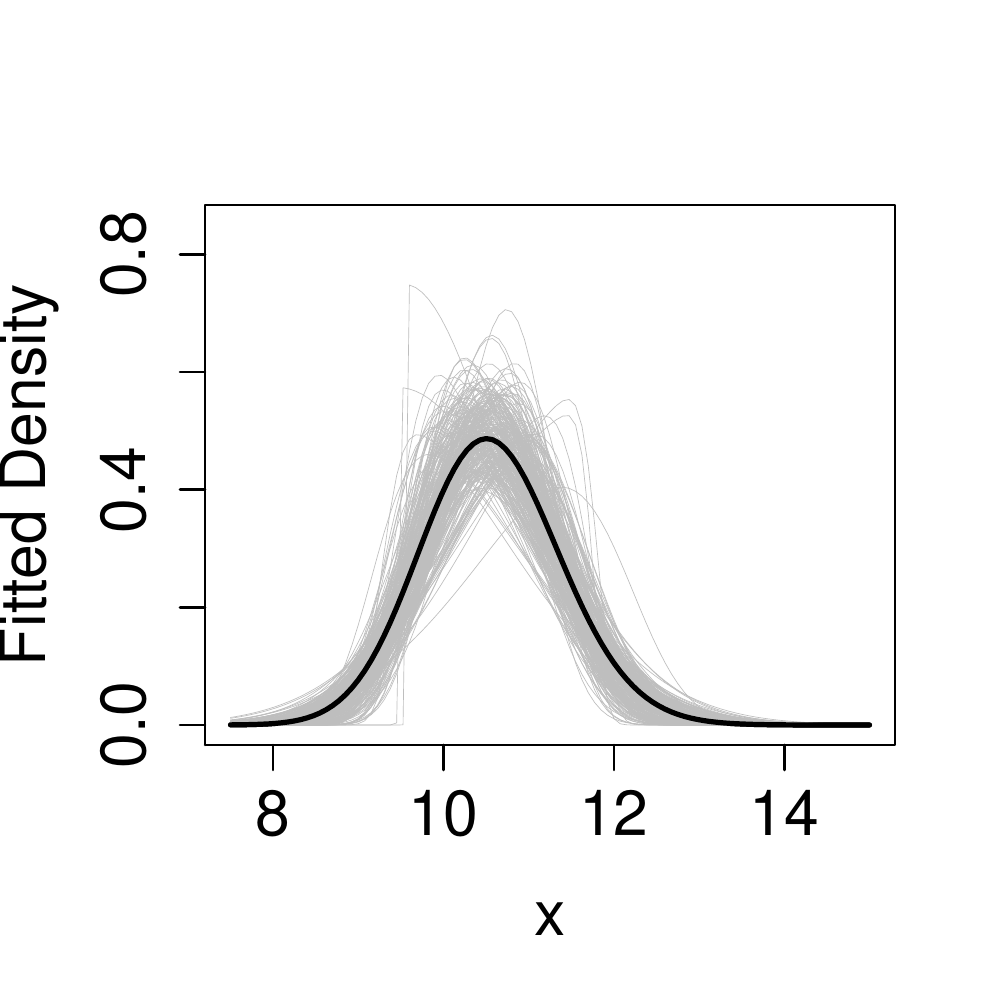} & 
\includegraphics[scale=0.5]{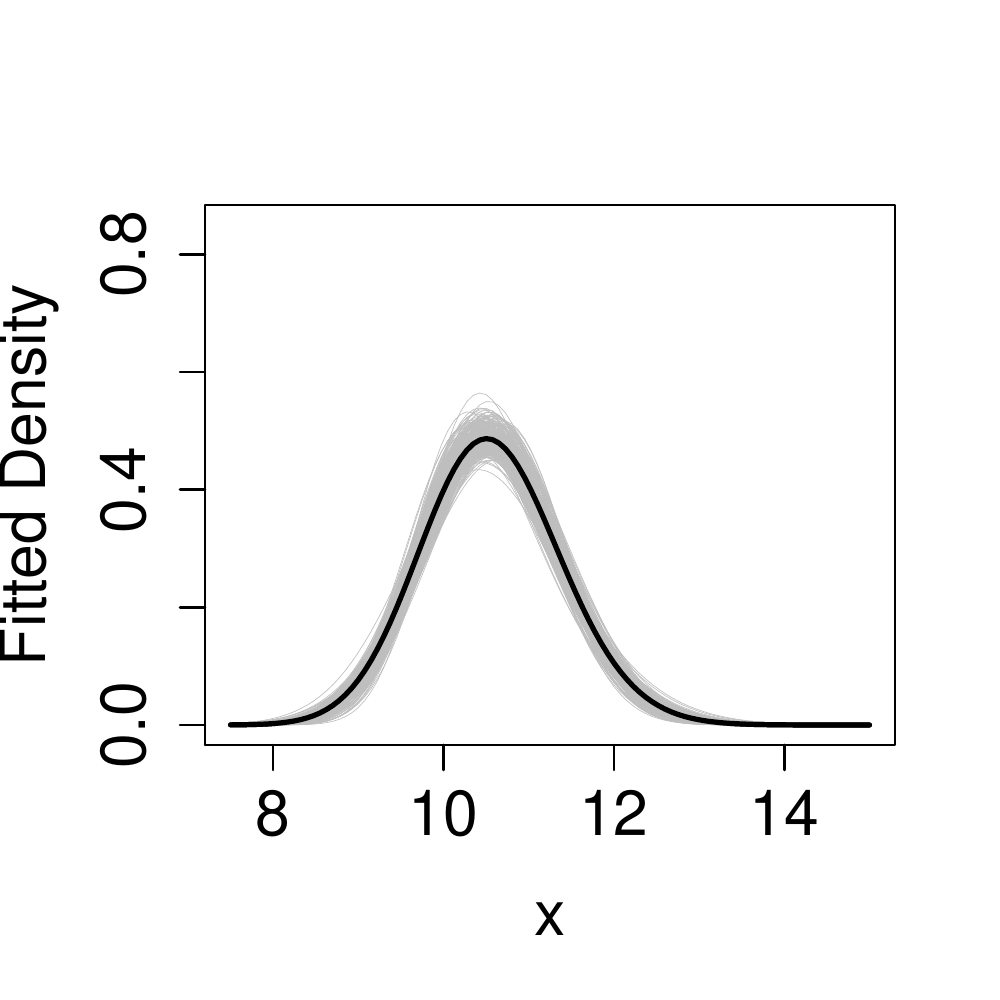} & 
\includegraphics[scale=0.5]{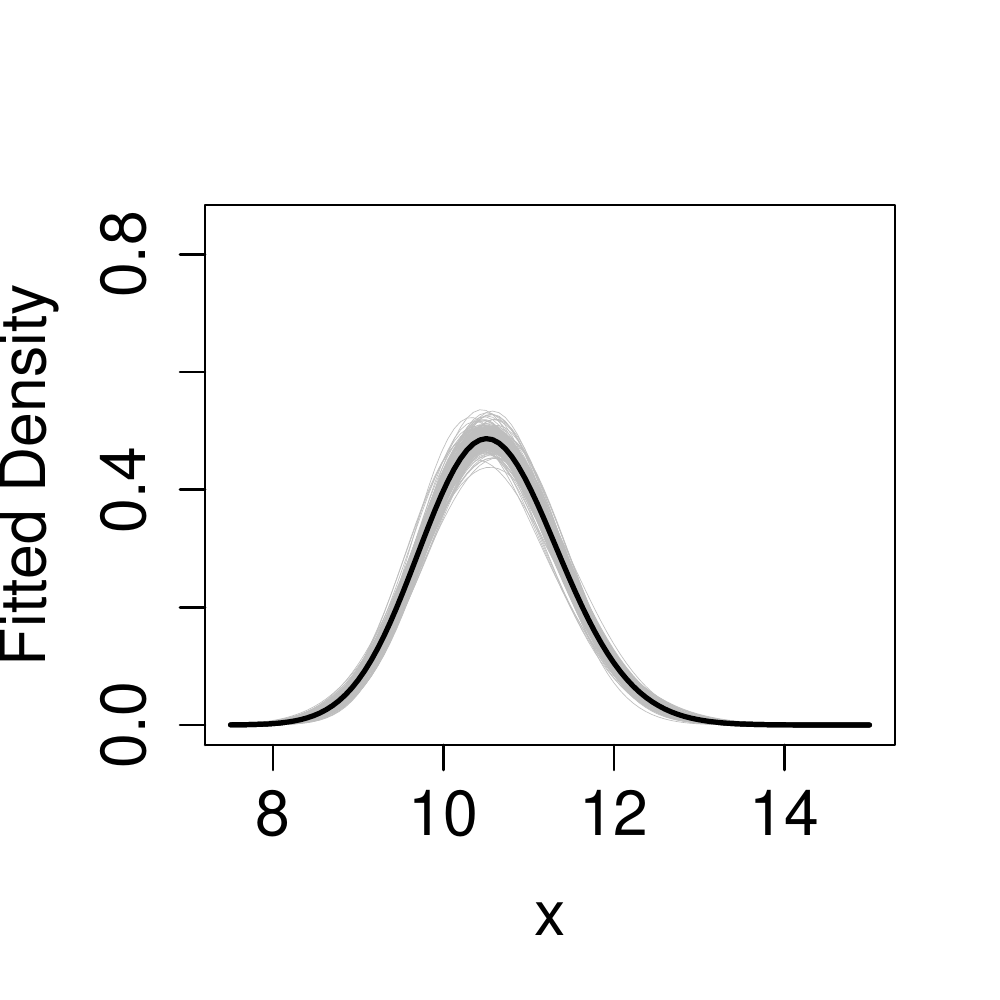} \\
(a) & (b) & (c) \\
\includegraphics[scale=0.5]{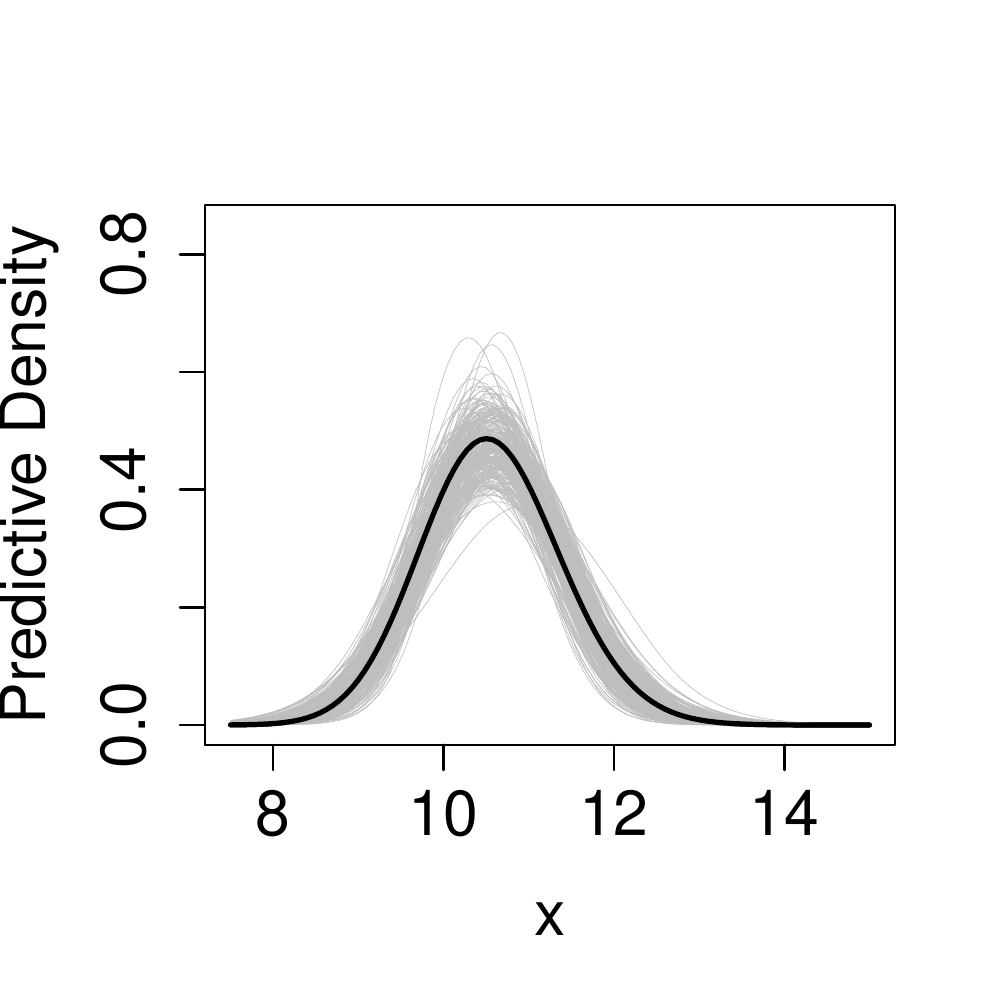} & 
\includegraphics[scale=0.5]{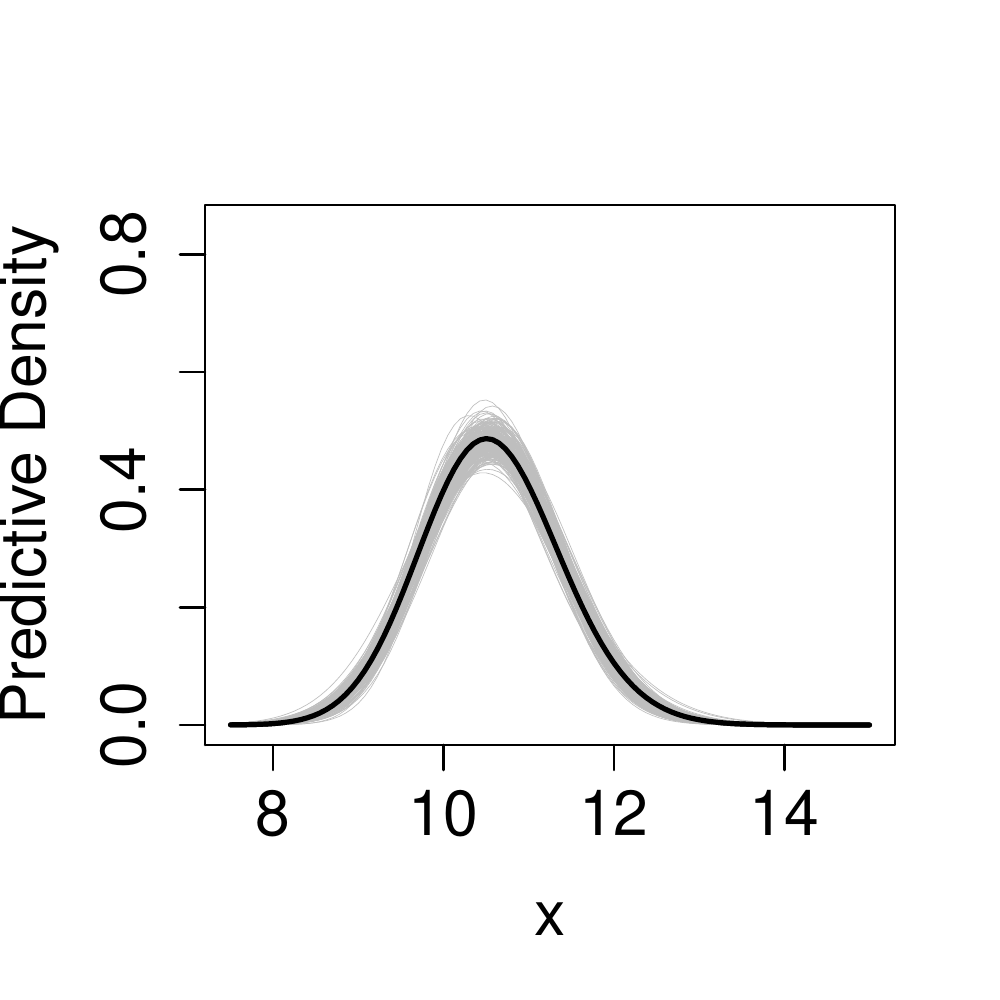} & 
\includegraphics[scale=0.5]{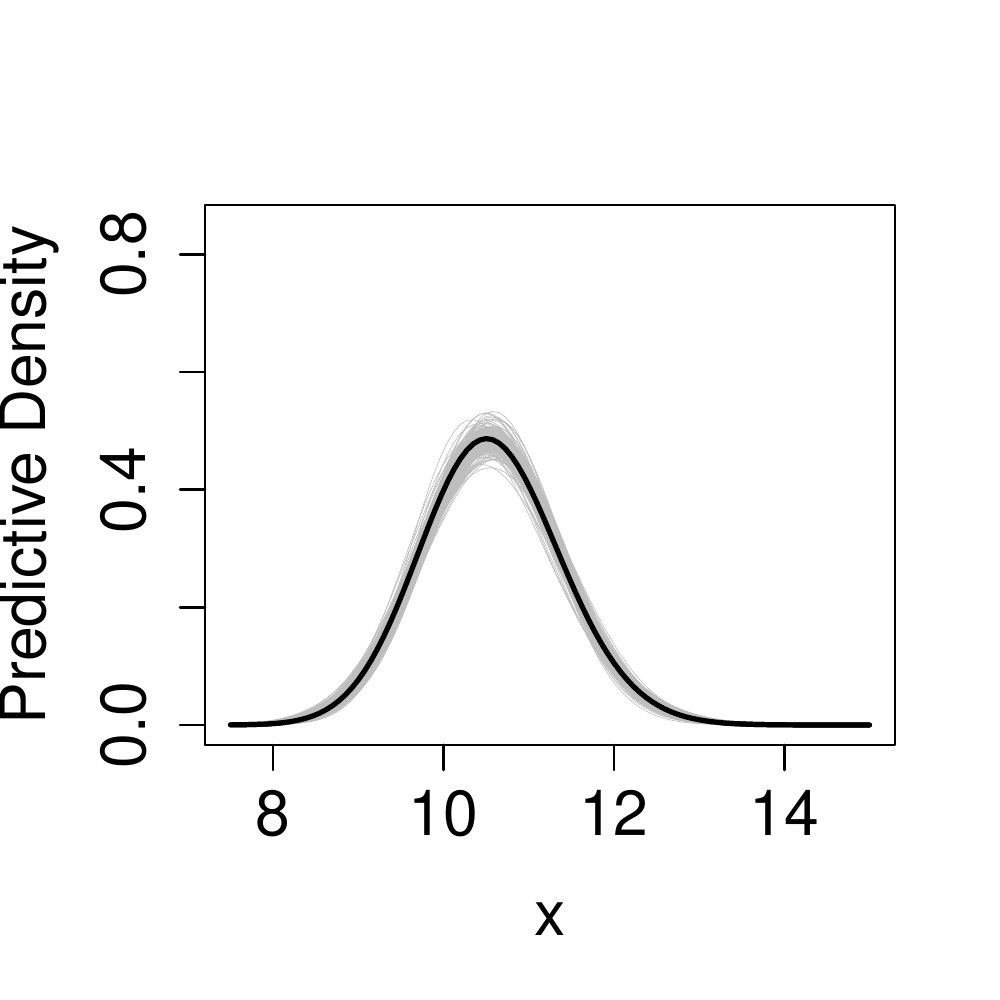} \\
(d) & (e) & (f) \\
\includegraphics[scale=0.5]{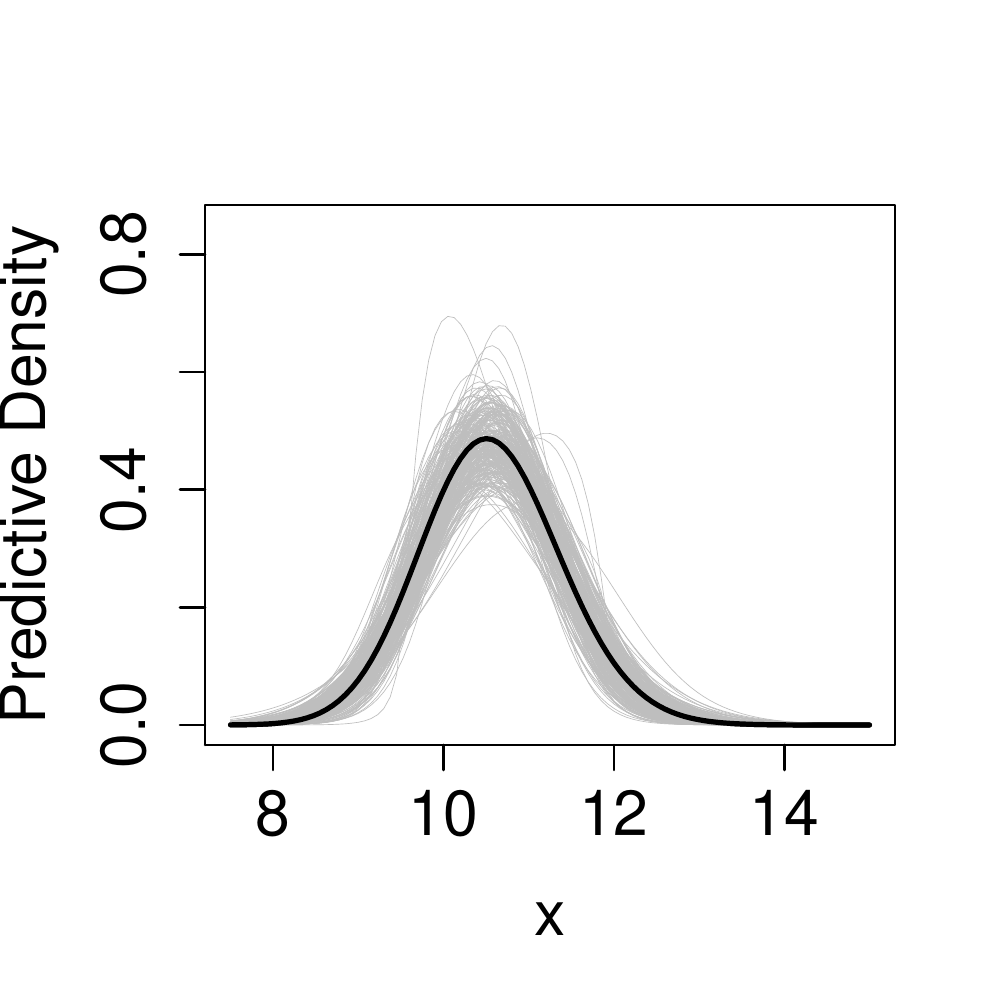} & 
\includegraphics[scale=0.5]{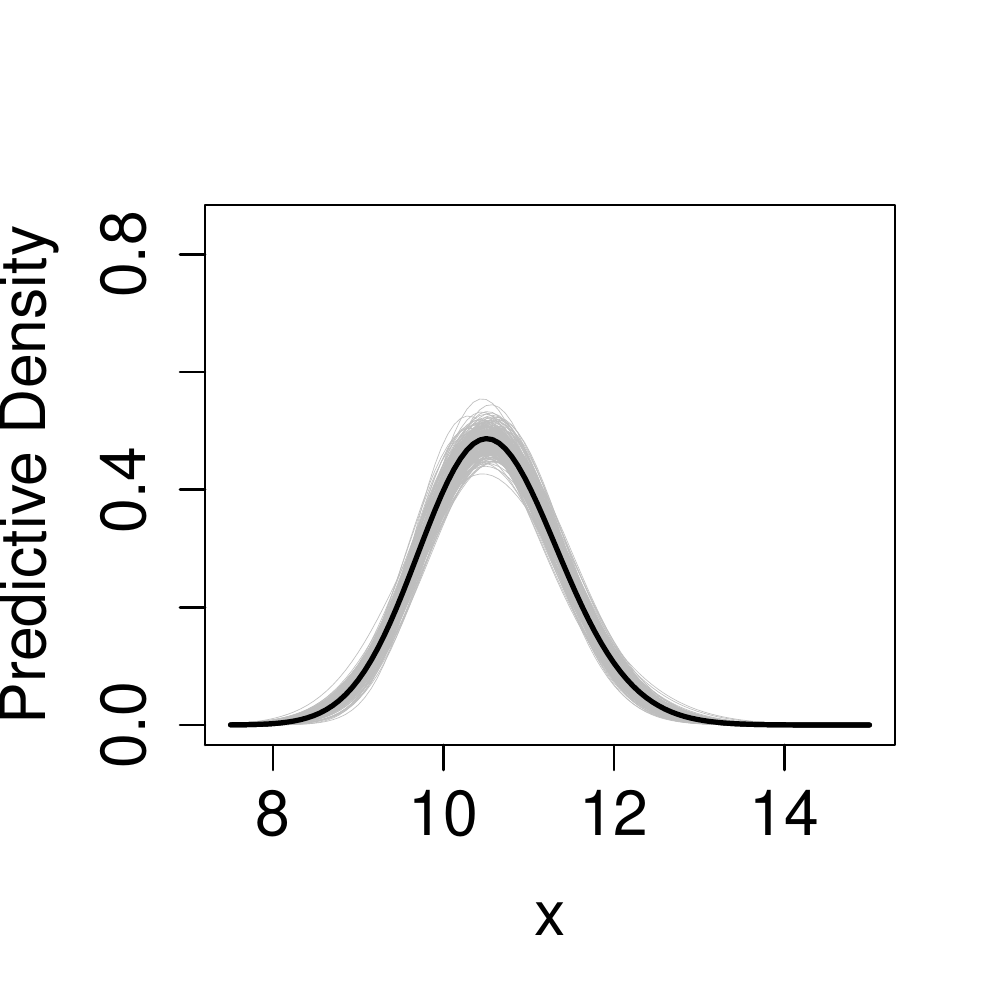} & 
\includegraphics[scale=0.5]{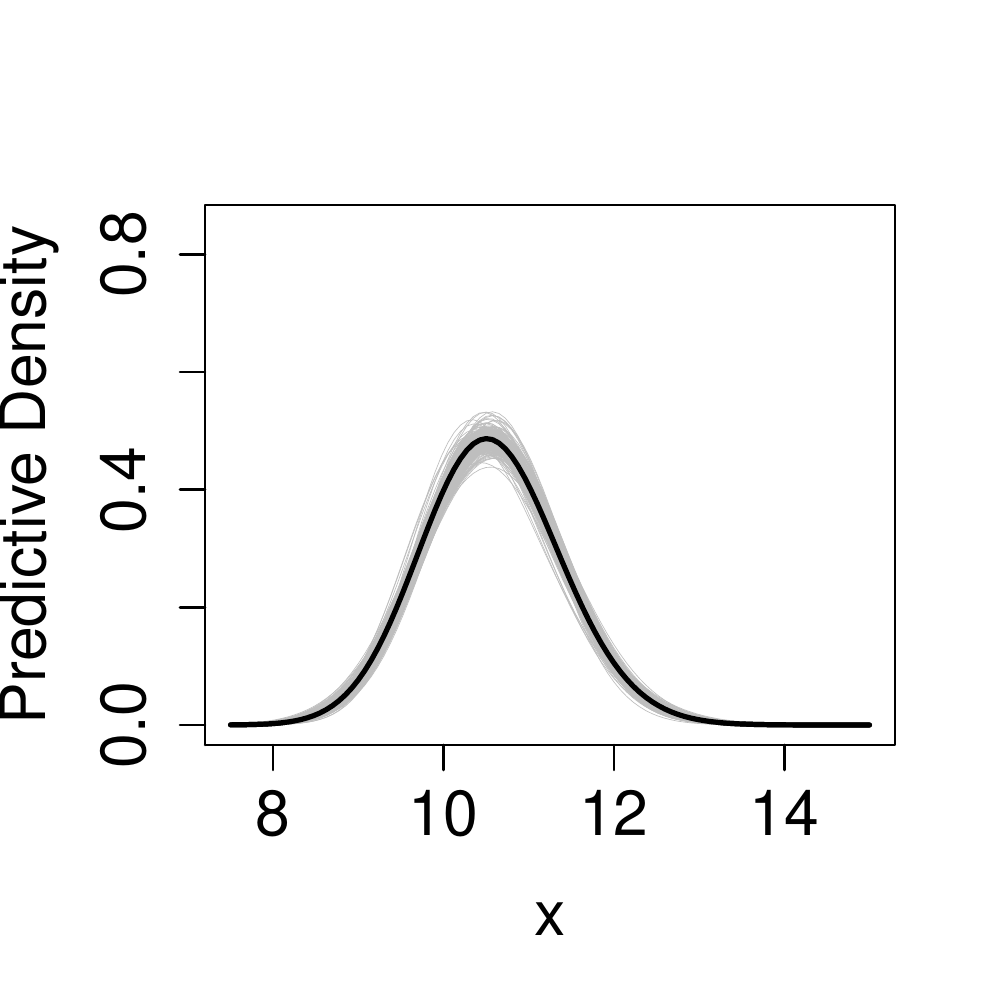} \\
(g) & (h) & (i) \\
\end{tabular}
\caption{(a) MLE-fitted densities ($\lambda = 1$, $n=50$), (b) MLE-fitted densities ($\lambda = 1$, $n=250$), (c) MLE-fitted densities ($\lambda = 1$, $n=500$); (d) Predictive densities for the independence Wasserstein prior ($\lambda = 1$, $n=50$), (e) Predictive densities for the independence Wasserstein prior ($\lambda = 1$, $n=250$), (f) Predictive densities for the independence Wasserstein prior ($\lambda = 1$, $n=500$); (g) Predictive densities for the independence Jeffreys prior ($\lambda = 1$, $n=50$), (d) Predictive densities for the independence Jeffreys prior ($\lambda = 1$, $n=250$), (d) Predictive densities for the independence Jeffreys prior ($\lambda = 1$, $n=500$). }
\label{fig:lambda1}
\end{figure}

Tables \ref{tab:lambda1}--\ref{tab:lambda5} in the Appendix show the results associated to the first simulation scenario. From Table \ref{tab:lambda1} in the Appendix, we observe that (in the case $\lambda = 1$) the estimation of the parameter $\lambda$ is indeed quite challenging for all sample sizes as the true model is very close to symmetry. 
Both priors (independence Jeffreys and independence Wasserstein) induce a marked shrinkage of the parameter $\lambda$ towards zero as the likelihood is relatively flat. This shrinkage naturally induces a bias in the Bayesian point estimators (posterior mean) for both priors. 
Although, for $n=50$, the coverage produced by the Jeffreys prior is slightly better than that produced by the Wasserstein prior, the average RMSE and standard deviations of the Bayes estimators associated to the Jeffreys prior are much larger.
This is likely a consequence of the very heavy tails of the Jeffreys prior which, together with the flatness of the likelihood, produce a heavy tailed posterior. Indeed, the MLE also exhibits a very large RMSE for $n=50$. 
The stronger regularisation induced by the Wasserstein prior also produces a faster concentration of the predictive posterior densities around the true model. The fit of the posterior predictive pdfs is particularly better than that obtained with the fitted pdfs using the MLEs (Figure \ref{fig:lambda1}). 
The cases $\lambda = 3,5$ (Tables \ref{tab:lambda3}--\ref{tab:lambda5} and Figures \ref{fig:lambda1}--\ref{fig:lambda5} in the Appendix) show that the estimation of the parameter $\lambda$ is much better behaved when the true value of $\lambda$ is away from $\lambda=0$, and the density function is clearly asymmetric. 
The performance of the independence Jeffreys and the independence Wasserstein in terms of all measures is quite similar. Since the true value of the parameter lies in the tails of the prior, the shrinkage effect of the priors is minimal.
In those cases, the MLE also exhibits a much larger RMSE for $n=50$.

Table \ref{tab:lrm} shows the results associated to the second simulation scenario. We notice that the performance of the Wasserstein prior is good for all sample sizes in terms of the chosen measures. Indeed, given that the prior is flat, the performance of the MLE coincides with that of the maximum a posteriori (MAP).
 
\section{Discussion}\label{sec:discussion}
We have introduced the Wasserstein prior, a prior based on the Wasserstein information matrix, which is invariant under reparameterisations. We have briefly discussed the link of the construction of this prior with concepts from information geometry.
We have also introduced the independence Wasserstein prior, which aims at reducing the functional dependence between the parameters in a similar fashion as the independence Jeffreys prior (and more generally, the reference prior \citep{yang:1997}).
The simulation study (results presented in the Appendix) shows that the Wasserstein prior induces a posterior with good frequentist properties (at least for the models studied here), compared to the posteriors induced by the Jeffreys prior and the fitted models using maximum likelihood estimation. 
Additional numerical examples related to the models presented here can be found at \url{https://github.com/FJRubio67/PIW}. 

As discussed in the introduction, objective priors are based on formal rules with specific aims.
The Wasserstein prior is based on a formal rule aiming at obtaining a prior that is invariant under reparameterisations. Consequently, the construction of such prior does not necessarily penalise model complexity, and thus may produce suboptimal results in sparse scenarios, such as linear regression models with many spurious variables. 

Natural extensions of our work include the calculation of the Wasserstein prior for other univariate continuous distributions (with bounded support, with positive support or supported on the entire real line). In this paper, we have taken a conservative position as we do not claim superiority of the Wasserstein prior over the Jeffreys prior in terms of a specific optimality criterion, even though the simulation study illustrates a competitive performance. Our work represents a step forward in the analysis of invariant priors obtained by a formal rule, and shows that it is possible to go beyond those induced by the Kullback-Leibler divergence. We believe it would be interesting to provide a theoretical treatment of the inferential properties of the Wasserstein prior, beyond the propriety of the posterior shown here. This includes the study of the asymptotic normality of the posterior distribution; establishing more formal links of the Wasserstein prior with information geometry \citep{kass:1996,kass:1989,nielsen:2020}; and the effect of the parameterisation on the orthogonality of parameters \citep{cox:1987} based on the Wasserstein information matrix. 

\clearpage
\section*{Appendix}

\section*{The Exponential Distribution}
Consider the exponential distribution with scale parameter $\theta >0$. The corresponding cdf and pdf are given by
\begin{eqnarray}\label{eq:exp}
F(x \mid \theta) &=& 1-\exp\left\{ -\dfrac{x}{\theta}\right\}, \quad f(x \mid \theta) = \dfrac{1}{\theta}\exp\left\{ -\dfrac{x}{\theta}\right\}, \quad x > 0.
\end{eqnarray}
In this case, the WIM and the Wasserstein prior are
\begin{eqnarray}\label{eq:piwexp}
W(\theta) = 2, \quad \pi_W(\theta) \propto 1.
\end{eqnarray}
Let $\bx = (x_1,\dots, x_n)^{\top}$ be an \textit{i.i.d.}~sample from \eqref{eq:exp}. Then, the posterior distribution of $\btheta$ associated to the Wasserstein prior \eqref{eq:piwexp} is proper if the sample size $n>1$. We omit the proof of this result, for the sake of space, as it is straightforward.

\section*{Proof of Theorem \ref{th:locsc}}
Let $X$ denote a random variable with cdf and pdf $F$ and $f$.
First, note that the first partial derivatives of the cdf are given by
\begin{eqnarray*}
\dfrac{\partial}{\partial \mu} F(x\mid\mu,\sigma) &=& - \dfrac{1}{\sigma}f_0\left( \dfrac{x-\mu}{\sigma} \right),\quad
 \dfrac{\partial}{\partial \sigma} F(x\mid\mu,\sigma) =  -\dfrac{x-\mu}{\sigma^2}f_0\left( \dfrac{x-\mu}{\sigma} \right).
\end{eqnarray*}
Then, the entries of the WIM are given by
\begin{eqnarray*}
W_{\mu,\mu} &=& \bbexp \left[  \dfrac{\left\{ \dfrac{\partial}{\partial \mu} F(X\mid\mu,\sigma) \right\}^2}{ \left\{ \dfrac{1}{\sigma}f_0\left( \dfrac{X-\mu}{\sigma}\right) \right\}^2} \right] \\
&=& 1.\\
W_{\sigma,\sigma} &=& \bbexp \left[  \dfrac{\left\{ \dfrac{\partial}{\partial \sigma} F(X\mid\mu,\sigma) \right\}^2}{ \left\{ \dfrac{1}{\sigma}f_0\left( \dfrac{X-\mu}{\sigma}\right) \right\}^2} \right] \\
&=& \bbexp\left[ \left\{\dfrac{X-\mu}{\sigma} \right\}^2\right] \quad  \left( \text{change of variable } t = \dfrac{x-\mu}{\sigma} \right)\\
&=& \int_{-\infty}^{\infty} t^2 f_0(t)dt.\\
W_{\mu,\sigma} &=& \bbexp \left[  \dfrac{\left\{ { \dfrac{\partial}{\partial \mu} F(X\mid\mu,\sigma) \dfrac{\partial}{\partial \sigma}} F(X\mid\mu,\sigma) \right\}}{ \left\{ \dfrac{1}{\sigma}f_0\left( \dfrac{X-\mu}{\sigma}\right) \right\}^2} \right] \\
 &=& \bbexp \left[ \dfrac{X-\mu}{\sigma} \right]\\
 &=& 0.
\end{eqnarray*}
Consequently,  the WIM and the Wasserstein prior for the location-scale family of $(\mu,\sigma)$ are given by
\begin{eqnarray*}
W(\mu,\sigma) = \left(
\begin{array}{c c}
1 & 0 \\
0 &  \int_{-\infty}^{\infty} t^2 f_0(t)dt.
\end{array}
\right),\quad \pi_W(\mu,\sigma) \propto 1.
\end{eqnarray*}

\section*{Proof of Theorem \ref{th:properls} }
Let $\blambda = (\lambda_1,\dots,\lambda_n)^{\top}$ and $\Lambda = \mbox{diag} (\lambda_1,\dots,\lambda_n)$. The posterior distribution is proper if the marginal likelihood (normalising constant) is finite. This is, we need to prove that
\begin{eqnarray*}
m(\bx) &=& \int_{\R \times \R_+^{n+1}} \prod_{i=1}^n \dfrac{\lambda_i^{1/2}}{\sqrt{2\pi}\sigma} \exp \left\{ -\dfrac{\lambda_i}{2\sigma^2}(x_i-\mu)^2 \right\} d\mu d\sigma dH(\blambda) \\
&=& \dfrac{1}{(2\pi)^{n/2}}\int_{\R \times \R_+^{n+1}}\left[ \prod_{i=1}^n \lambda_i^{1/2} \right]\dfrac{1}{\sigma^n} \exp \left\{ -\dfrac{1}{2\sigma^2}\sum_{i=1}^n \lambda_i(x_i-\mu)^2 \right\} d\mu d\sigma dH(\blambda)< \infty.
\end{eqnarray*}
Notice first that
\begin{eqnarray*}
\sum_{i=1}^n \lambda_i(x_i-\mu)^2  = A (\mu - B)^2 + s(\blambda),
\end{eqnarray*}
where $A = \sum_{i=1}^n \lambda_i$, $B= \dfrac{\sum_{i=1}^n \lambda_i x_i}{\sum_{i=1}^n \lambda_i}$,  $s(\blambda) = \bx^{\top} \Lambda \bx - \bx^{\top} \Lambda \bone (\bone^{\top} \Lambda \bone)^{-1} \bone^{\top} \Lambda \bx$, and $\bone$ is a vector of $1$s of length $n$. Using this decomposition and integrating out $\mu$ as a normal distribution, we obtain
\begin{eqnarray*}
m(\bx) &=& \dfrac{1}{(2\pi)^{\frac{n-1}{2}}}\int_{\R_+^{n}}   \left[ \prod_{i=1}^n \lambda_i^{1/2} \right ] (\bone^{\top}\Lambda \bone)^{-1/2} \dfrac{1}{\sigma^{n-1}} \exp \left\{ -\dfrac{s(\blambda)}{2\sigma^2} \right\} d\sigma dH(\blambda).
\end{eqnarray*}
Now, for $n\geq3$, integrating this expression with respect to $\sigma^2$ as a Gamma distribution we obtain,
\begin{eqnarray*}
m(\bx) &=& \dfrac{1}{(2\pi)^{\frac{n-1}{2}}}\int_{\R_+^{n+1}}   \left[ \prod_{i=1}^n \lambda_i^{1/2} \right ] (\bone^{\top}\Lambda \bone)^{-1/2} 2^{\frac{n}{2}-2} \Gamma \left( \dfrac{n-2}{2} \right) s(\blambda)^{-\frac{n-2}{2}}  dH(\blambda).
\end{eqnarray*}
Now, using Lemma 3 from \cite{fernandez:2000}, $s(\blambda)$ has a lower bound proportional to $\lambda_{(n)} = \max_{i=1,\dots,n} \{\lambda_1,\dots,\lambda_n\}$. Moreover, by Lemma 1 from \cite{fernandez:2000}, $\bone^{\top}\Lambda \bone$ has a lower bound proportional to $\lambda_{(n)}$. Consequently, there exists a constant $K$ such that
\begin{eqnarray*}
\left[\prod_{i=1}^n\lambda_i^{1/2}\right](\bone^{\top}\Lambda \bone)^{-1/2} s(\blambda)^{-\frac{n-2}{2}}  \leq K\lambda_{(1)}^{1/2} = K\min\{ \lambda_1,\dots,\lambda_n\}^{1/2}\leq K\lambda_i^{1/2} \,\,\, \text{for all }i = 1,\dots,n.
\end{eqnarray*}
Consequently, there exists a constant $0 < \tilde{K} < \infty$ such that,
\begin{eqnarray*}
m(\bx) &\leq& \tilde{K}\int_{\R_+}    \lambda^{1/2} dH(\lambda) < \infty.
\end{eqnarray*}

\section*{Proof of Theorem \ref{th:sn}}
\begin{itemize}

\item[(i)] The skew normal cdf can be written as \citep{azzalini:1985}  $S(x \mid \alpha) = \Phi(x) - 2 T(x,\alpha)$, where $T(x,\alpha)$ is the Owen's T function 
$$T(x,a)={\frac {1}{2\pi }}\int _{0}^{\alpha}{\frac {e^{-{\frac {1}{2}}x^{2}(1+t^{2})}}{1+t^{2}}}dt\quad \left(-\infty < x, \alpha <+\infty \right).$$
Using the Fundamental Theorem of Calculus
\begin{eqnarray*}
\dfrac{\partial}{\partial \alpha}S(x \mid \alpha) &=& - 2\dfrac{\partial}{\partial \alpha}T(x,\alpha) = - \frac{1}{\pi} \frac {e^{-{\frac {1}{2}}x^{2}(1+\alpha^{2})}}{1+\alpha^{2}}.
\end{eqnarray*}
Replacing this expression in formula \eqref{eq:WIM} together with the relationship $\Phi (x)={\frac  {1}{2}}\left[1+\operatorname {erf}\left({\frac  {x}{{\sqrt  {2}}}}\right)\right]$ we obtain
\begin{eqnarray*}
W(\alpha) = \int_{-\infty}^{\infty} \frac{\sqrt{2} e^{-\frac{1}{2} \left(2 \alpha ^2+1\right) x^2}}{\pi ^{3/2} \left(\alpha
   ^2+1\right)^2 \left(\mbox{erf}\left(\frac{\alpha  x}{\sqrt{2}}\right)+1\right)}dx.
\end{eqnarray*}

\item[(ii)] The symmetry of $\pi_W$ is a consequence of the integrand of $W$ being a function of $\alpha^2$ and the $\mbox{erf}$ function applied to $\alpha$, together with the property $\mbox{erf}(z) = - \mbox{erf}(-z)$. 

\item[(iii)]  Note that for $\alpha=0$
\begin{eqnarray*}
W(0) = \int_{-\infty}^{\infty} \frac{\sqrt{2} e^{-\frac{1}{2} x^2}}{\pi ^{3/2} }dx < \infty.
\end{eqnarray*}
Now, for $M>0$ and $\alpha \in [-M,M]$, and using the symmetry of $W$, there exists $K_1>0$ such that
\begin{eqnarray*}
W(\alpha) &=& \int_{-\infty}^{\infty} \frac{\sqrt{2} e^{-\frac{1}{2} \left(2 \alpha ^2+1\right) x^2}}{\pi ^{3/2} \left(\alpha
   ^2+1\right)^2 \left(\mbox{erf}\left(\frac{\alpha  x}{\sqrt{2}}\right)+1\right)}dx \\
   &=&\dfrac{K_1 \sqrt{2}}{\pi ^{3/2} \left(\alpha^2+1\right)^2}  \int_{-\infty}^{\infty} \dfrac{\phi(2\alpha x)^2 \phi(x)}{\Phi(\alpha x)} dx\\
      &=& \dfrac{2K_1 \sqrt{2}}{ \pi ^{3/2}\left(\alpha^2+1\right)^2}  \int_{0}^{\infty} r(\alpha x)\phi(2\alpha x) \phi(x) dx,
\end{eqnarray*}
where $r(x) = \dfrac{\phi(x)}{\Phi(x)} $ denotes the inverse Mills ratio. It is well known that $r(x) $ is a decreasing function and that $r(x) \sim \vert x \vert$ as $x \to -\infty$ and $r(x) \sim \vert 1/x \vert$ as $x \to \infty$. Then, for $\alpha \in [-M,M]$ and $x>0$, there exists $K_2>0$ such that $r(\alpha x)\phi(2\alpha x) \leq K_2$ and
\begin{eqnarray*}
W(\alpha) &\leq& \dfrac{K_1 K_2 \sqrt{2}}{ \pi ^{3/2}\left(\alpha^2+1\right)^2} \leq \dfrac{K_1 K_2 \sqrt{2}}{ \pi ^{3/2}}.
\end{eqnarray*}
Since $W(\alpha)$ is upper-bounded by a finite constant for $\vert \alpha \vert \leq M$, it follows that
\begin{eqnarray*}
\int_{-M}^M \pi_W(\alpha) d\alpha \propto 2\int_{0}^M \sqrt{W(\alpha)} d\alpha< \infty.
\end{eqnarray*}
Now, for $ \alpha  > M$, and since $r(\cdot)$ and $\phi(\cdot)$ are decreasing functions, there exists $K_3 >0$ such that $r(\alpha x)\phi(2\alpha x) \leq K_3$, and
\begin{eqnarray*}
W(\alpha)    &\leq& \dfrac{K_1 K_3 \sqrt{2}}{ \pi ^{3/2}\left(\alpha^2+1\right)^2}.
\end{eqnarray*}
Consequently,
\begin{eqnarray*}
\int_{M}^{\infty} \pi_W(\alpha) d\alpha &\propto& \int_{M}^{\infty} \sqrt{W(\alpha)} d\alpha \leq \int_{M}^{\infty} \sqrt{ \dfrac{K_1 K_3 \sqrt{2}}{ \pi ^{3/2}\left(\alpha^2+1\right)^2}}  d\alpha< \infty.
\end{eqnarray*}
Finally, for $\alpha < - M$, appealing to the symmetry of $\pi_W$, we also obtain that
\begin{eqnarray*}
\int_{-\infty}^{-M} \pi_W(\alpha) d\alpha < \infty. 
\end{eqnarray*}

\item[(iv)]  Using the expression
\begin{eqnarray*}
W(\alpha) &=&\dfrac{2K_1 \sqrt{2}}{ \pi ^{3/2}\left(\alpha^2+1\right)^2}  \int_{0}^{\infty} r(\alpha x)\phi(2\alpha x) \phi(x) dx,
\end{eqnarray*}
for $\alpha>0$ and noting that $\phi(x)$ is upper bounded, it follows that
\begin{eqnarray*}
W(\alpha) &\leq&\dfrac{2K_1 \sqrt{2}}{ \pi ^{3/2}\left(\alpha^2+1\right)^2}  \int_{0}^{\infty} r(\alpha x)\phi(2\alpha x)  dx.
\end{eqnarray*}
Consider now the change of variable $t=\alpha x$. Then, we obtain the upper bound
\begin{eqnarray}\label{eq:upb}
W(\alpha) &\leq&\dfrac{2K_1 \sqrt{2}}{ \pi ^{3/2}\left(\alpha^2+1\right)^2\alpha}  \int_{0}^{\infty} r(t)\phi(2t)  dt.
\end{eqnarray}
Now consider the change of variable $t=\alpha x$ applied to the expression
\begin{eqnarray*}
W(\alpha) &=&\dfrac{2K_1 \sqrt{2}}{ \pi ^{3/2}\left(\alpha^2+1\right)^2}  \int_{0}^{\infty} r(\alpha x)\phi(2\alpha x) \phi(x) dx\\
 &=&\dfrac{2K_1 \sqrt{2}}{ \pi ^{3/2}\left(\alpha^2+1\right)^2\alpha}  \int_{0}^{\infty} r(t)\phi(2t) \phi\left(\frac{t}{\alpha}\right) dt.
\end{eqnarray*}
Note now that for $\alpha \geq A > 0$, $ \phi\left(\frac{t}{\alpha}\right) \geq \phi\left(\frac{t}{A}\right)$, for all $t>0$. Then,
\begin{eqnarray}\label{eq:lowb}
W(\alpha) &\geq& \dfrac{2K_1 \sqrt{2}}{ \pi ^{3/2}\left(\alpha^2+1\right)^2\alpha}  \int_{0}^{\infty} r(t)\phi(2t) \phi\left(\frac{t}{A}\right) dt.
\end{eqnarray}
Combining \eqref{eq:upb}-\eqref{eq:lowb} we obtain that $W(\alpha)$ has tails of order ${\mathcal O}(\vert \alpha \vert^{-5})$. This implies that $\pi_W(\alpha)$ has tails of order ${\mathcal O}(\vert \alpha \vert^{-5/2})$.

\end{itemize}

\section*{Proof of Theorem \ref{th:propersn}}
The marginal likelihood can be upper bounded by
\begin{eqnarray*}
m(\bx) &=& \int_{{\mathbb R}^2 \times {\mathbb R}_+} \prod_{i=1}^n \dfrac{2}{\sigma} \phi\left(\dfrac{x_i-\mu}{\sigma}\right) \Phi\left(\alpha \dfrac{x_i-\mu}{\sigma}\right) \pi(\mu,\sigma,\alpha) d\mu d\alpha d\sigma\\
&\leq &\int_{{\mathbb R}^2 \times {\mathbb R}_+} \prod_{i=1}^n \dfrac{2}{\sigma} \phi\left(\dfrac{x_i-\mu}{\sigma}\right) \pi(\mu,\sigma,\alpha) d\mu d\alpha d\sigma \\
(\text{integrating out } \alpha) &\propto &\int_{{\mathbb R} \times {\mathbb R}_+} \prod_{i=1}^n \dfrac{2}{\sigma} \phi\left(\dfrac{x_i-\mu}{\sigma}\right) \pi(\mu,\sigma) d\mu d\sigma.
\end{eqnarray*}
The last expression is proportional to the marginal likelihood associated to a normal sampling model together with the Wasserstein prior. By Theorem \ref{th:properls}, we have that this marginal likelihood is finite.

\section*{Proof of Theorem \ref{th:lrm}}
First, note that the cdf and pdf associated to the $i$th observation of the normal linear regression model are given by
\begin{eqnarray*}
F(y_i \mid \bx_i, \bbeta,\sigma) &=& \Phi\left( \dfrac{y_i - \bx_i^{\top}\bbeta }{\sigma} \right), \quad
f(y_i \mid \bx_i, \bbeta,\sigma) = \dfrac{1}{\sigma}\phi\left( \dfrac{y_i - \bx_i^{\top}\bbeta }{\sigma} \right), 
\end{eqnarray*}
where $\Phi$ and $\phi$ are the standard normal cdf and pdf, respectively. The derivatives with respect to the parameters are 
\begin{eqnarray*}
\dfrac{\partial}{\partial \bbeta_k}F(y_i \mid \bx_i, \bbeta,\sigma) &=& - \dfrac{\bx_{ik}}{\sigma} \phi\left( \dfrac{y_i - \bx_i^{\top}\bbeta }{\sigma} \right), \,\,\,\,\, k=1,\dots, p, \\
\dfrac{\partial}{\partial \sigma}F(y_i \mid \bx_i, \bbeta,\sigma) &=&- \dfrac{y_i-\bx_i^{\top}\bbeta}{\sigma^2}\phi\left( \dfrac{y_i - \bx_i^{\top}\bbeta }{\sigma} \right). 
\end{eqnarray*}
Replacing these expressions in formula \eqref{eq:WIM}, we obtain
\begin{eqnarray*}
I_{\bbeta_j \bbeta_k} (\bbeta,\sigma) &=& \bx_{ij}\bx_{ik}, \\
I_{\bbeta_j \sigma} (\bbeta,\sigma) &=& 0,\\
I_{\sigma,\sigma} (\bbeta,\sigma) &=& \int_0^\infty u^2 \phi(u)du = 1.
\end{eqnarray*}
Consequently, using Proposition 5 in \cite{li:2019} and the assumption of independence of the errors, the WIM for the entire sample is given by
\begin{eqnarray*}
W(\bbeta,\sigma) &=& \left(
\begin{array}{c c}
\bX^{\top}\bX & 0 \\
0 &  1
\end{array}
\right).
\end{eqnarray*}
Taking the square root of the determinant of the WIM, we obtain the Wasserstein prior $\pi_W(\bbeta,\sigma) \propto 1$.

\section*{Proof of Theorem \ref{th:properlrm}}
The posterior distribution is proper if the marginal likelihood (normalising constant) is finite. This is, we need to prove that
\begin{eqnarray*}
m(\by \mid \bx) &=& \int_{\R^p \times \R_+} \prod_{i=1}^n \dfrac{1}{\sqrt{2\pi}\sigma} \exp \left\{ -\dfrac{1}{2\sigma^2}(y_i - \bx_i^{\top}\bbeta)^2 \right\} d\bbeta d\sigma  \\
&=& \dfrac{1}{(2\pi)^{n/2}}\int_{\R^p \times \R_+} \dfrac{1}{\sigma^n} \exp \left\{ -\dfrac{1}{2\sigma^2}(\by - \bX\bbeta)^{\top}(\by - \bX\bbeta) \right\} d\bbeta d\sigma< \infty.
\end{eqnarray*}
Consider the classical decomposition
\begin{eqnarray*}
(\by - \bX\bbeta)^{\top}(\by - \bX\bbeta) = (\by - \bX\widehat{\bbeta})^{\top}(\by - \bX\widehat{\bbeta}) + (\bbeta - \widehat{\bbeta})^{\top}(\bX^{\top}\bX)(\bbeta - \widehat{\bbeta}),
\end{eqnarray*}
where $\widehat{\bbeta} = (\bX^{\top}\bX)^{-1}\bX \by$. Replacing this expression in the marginal likelihood and integrating $\bbeta$ out as a $p$-variate normal distribution and $\sigma^2$ as a Gamma distribution, and using that $\bX$ has full column rank and that $\by $ is not in the column space of $\bX$, we obtain for $n>  p +1$ 
\begin{eqnarray*}
m(\bx) \leq K \mbox{det}(\bX^{\top}\bX)^{-\frac{1}{2}} [(\by - \bX\widehat{\bbeta})^{\top}(\by - \bX\widehat{\bbeta}) ]^{-\frac{n-p-1}{2}} < \infty,
\end{eqnarray*}
for a positive constant $K>0$.

\clearpage
\section*{Simulation Results}
Throughout, we denote by `mMean' the average of the posterior means across the $N$ simulated samples. Similarly, `mSD' represent the average of the posterior standard deviations; `mRMSE' denotes the average of the root mean squared errors; `Coverage' is the coverage proportion of the $95\%$ credible intervals; `mMLE' is the average of the maximum likelihood estimators; and `RMSE-MLE' is the root mean squared error of the maximum likelihood estimators across the $N$ samples.

\begin{table}[ht]
\centering
\begin{tabular}{ccccccc}
  \hline
 & $\mu$ (10) & $\sigma$ (1) & $\lambda$ (1)  & $\mu$ (10) & $\sigma$ (1) & $\lambda$ (1) \\ 
  \hline
  \multicolumn{7}{c}{  $n=50$} \\
  & \multicolumn{3}{c}{ Wasserstein} & \multicolumn{3}{c}{ Jeffreys} \\
  mMean & 10.502 & 0.983 & 0.218 & 10.487 & 1.019 & 0.303\\ 
  mSD & 0.474 & 0.179 & 1.254 & 0.509 & 0.194 & 2.510 \\   
  mRMSE & 0.726 & 0.213 & 1.691 & 0.766 & 0.231 & 3.045\\ 
  Coverage & 0.932 & 0.992 & 0.960 & 0.948 & 0.988 & 0.952 \\ 
  mMLE & 10.388 & 1.086 & 1.917 & -- & -- & -- \\ 
  RSME-MLE & 0.819 & 0.214 & 16.482 & -- & -- & --\\  
      \multicolumn{7}{c}{  $n=250$} \\
  mMean & 10.391 & 0.941 & 0.355 & 10.359 & 0.958 & 0.432\\ 
   mSD & 0.350 & 0.102 & 0.662 & 0.358 & 0.107 & 0.693 \\  
  mRMSE & 0.560 & 0.139 & 1.037 & 0.553 & 0.140 & 1.051 \\ 
  Coverage & 0.860 & 0.960 & 0.832 & 0.876 & 0.972 & 0.868\\ 
  mMLE & 10.163 & 1.010 & 0.793 & -- & -- & -- \\ 
  RSME-MLE & 0.464 & 0.101 & 0.887 & -- & -- & -- \\ 
        \multicolumn{7}{c}{  $n=500$} \\
      mMean & 10.346 & 0.932 & 0.419  & 10.297 & 0.945 & 0.512 \\ 
  mSD & 0.304 & 0.085 & 0.554 & 0.305 & 0.087 & 0.564 \\ 
  mRMSE & 0.493 & 0.127 & 0.895  & 0.463 & 0.124 & 0.850 \\ 
  Coverage & 0.808 & 0.924 & 0.804 & 0.876 & 0.932 & 0.880  \\ 
  mMLE & 10.099 & 0.994 & 0.866 & -- & -- & --  \\ 
  RSME-MLE & 0.329 & 0.081 & 0.605 & -- & -- & --  \\
   \hline
\end{tabular}
\caption{Simulation results for the independence Wasserstein prior and the independence Jeffreys prior.}
\label{tab:lambda1}
\end{table}

\begin{table}[ht]
\centering
\begin{tabular}{ccccccc}
  \hline
 & $\mu$ (10) & $\sigma$ (1) & $\lambda$ (3)  & $\mu$ (10) & $\sigma$ (1) & $\lambda$ (3) \\ 
  \hline
  \multicolumn{7}{c}{  $n=50$} \\
    & \multicolumn{3}{c}{ Wasserstein} & \multicolumn{3}{c}{ Jeffreys} \\
mMean & 10.401 & 0.843 & 1.976 & 10.324 & 0.881 & 3.359 \\ 
  mSD & 0.342 & 0.153 & 2.501 & 0.334 & 0.160 & 6.295 \\ 
  mRMSE & 0.560 & 0.258 & 3.537 & 0.508 & 0.248 & 7.273\\ 
  Coverage & 0.780 & 0.820 & 0.764 & 0.832 & 0.860 & 0.860 \\ 
  mMLE & 10.078 & 0.968 & 10.494 & -- & -- & --   \\ 
  RSME-MLE & 0.324 & 0.178 & 34.451 & -- & -- & --   \\
      \multicolumn{7}{c}{  $n=250$} \\
mMean & 10.046 & 0.972 & 2.884  & 10.032 & 0.981 & 3.029\\ 
  mSD & 0.092 & 0.073 & 0.756 & 0.085 & 0.071 & 0.776 \\ 
  mRMSE & 0.125 & 0.103 & 1.065 & 0.115 & 0.099 & 1.081 \\ 
  Coverage & 0.932 & 0.916 & 0.912  & 0.948 & 0.920 & 0.928\\ 
  mMLE & 10.007 & 0.994 & 3.138 & -- & -- & --   \\ 
  RSME-MLE & 0.069 & 0.067 & 0.827 & -- & -- & --   \\ 
        \multicolumn{7}{c}{  $n=500$} \\
mMean & 10.017 & 0.988 & 2.939 & 10.012 & 0.991 & 2.997 \\ 
   mSD & 0.049 & 0.048 & 0.484 & 0.049 & 0.048 & 0.495 \\ 
  mRMSE & 0.068 & 0.066 & 0.668 & 0.067 & 0.065 & 0.675 \\ 
  Coverage & 0.924 & 0.936 & 0.936& 0.920 & 0.940 & 0.940 \\ 
  mMLE & 10.005 & 0.996 & 3.041 & -- & -- & -- \\ 
  RSME-MLE & 0.049 & 0.047 & 0.510 & -- & -- & -- \\ 
   \hline
\end{tabular}
\caption{Simulation results for the independence Wasserstein prior and the independence Jeffreys prior.}
\label{tab:lambda3}
\end{table}

\begin{table}[ht]
\centering
\begin{tabular}{ccccccc}
  \hline
 & $\mu$ (10) & $\sigma$ (1) & $\lambda$ (5)  & $\mu$ (10) & $\sigma$ (1) & $\lambda$ (5) \\ 
  \hline
  \multicolumn{7}{c}{  $n=50$} \\
    & \multicolumn{3}{c}{ Wasserstein} & \multicolumn{3}{c}{ Jeffreys} \\
mMean & 10.272 & 0.865 & 4.702 & 10.232 & 0.880 & 5.306 \\ 
mSD & 0.269 & 0.144 & 5.596  & 0.271 & 0.150 & 8.338 \\
  mRMSE & 0.414 & 0.241 & 7.120  & 0.383 & 0.230 & 9.492\\ 
  Coverage & 0.804 & 0.824 & 0.816 & 0.876 & 0.860 & 0.892 \\ 
  mMLE & 10.037 & 0.968 & 35.886  & -- & -- & -- \\ 
  RSME-MLE & 0.199 & 0.161 & 73.255  & -- & -- & -- \\ 
      \multicolumn{7}{c}{  $n=250$} \\
mMean & 10.018 & 0.987 & 5.054  & 10.011 & 0.992 & 5.367\\ 
mSD & 0.049 & 0.057 & 1.318 & 0.048 & 0.058 & 1.479 \\ 
  mRMSE & 0.069 & 0.080 & 1.825 & 0.067 & 0.079 & 2.018 \\ 
  Coverage & 0.928 & 0.920 & 0.916 & 0.952 & 0.932 & 0.916 \\ 
  mMLE & 10.004 & 0.996 & 5.429  & -- & -- & -- \\ 
  RSME-MLE & 0.050 & 0.059 & 1.746  & -- & -- & -- \\ 
        \multicolumn{7}{c}{  $n=500$} \\
mMean & 10.009 & 0.993 & 4.973 & 10.005 & 0.996 & 5.112 \\ 
  mSD & 0.033 & 0.040 & 0.860 & 0.033 & 0.040 & 0.882 \\ 
  mRMSE & 0.047 & 0.055 & 1.174 & 0.046 & 0.055 & 1.209 \\ 
  Coverage & 0.912 & 0.948 & 0.928 & 0.928 & 0.936 & 0.924 \\ 
  mMLE & 10.002 & 0.998 & 5.142  & -- & -- & --\\ 
  RSME-MLE & 0.036 & 0.041 & 0.957 & -- & -- & -- \\
   \hline
\end{tabular}
\caption{Simulation results for the independence Wasserstein prior and the independence Jeffreys prior.}
\label{tab:lambda5}
\end{table}

\pagebreak 

\begin{figure}
\begin{tabular}{ccc}
\centering
\includegraphics[scale=0.5]{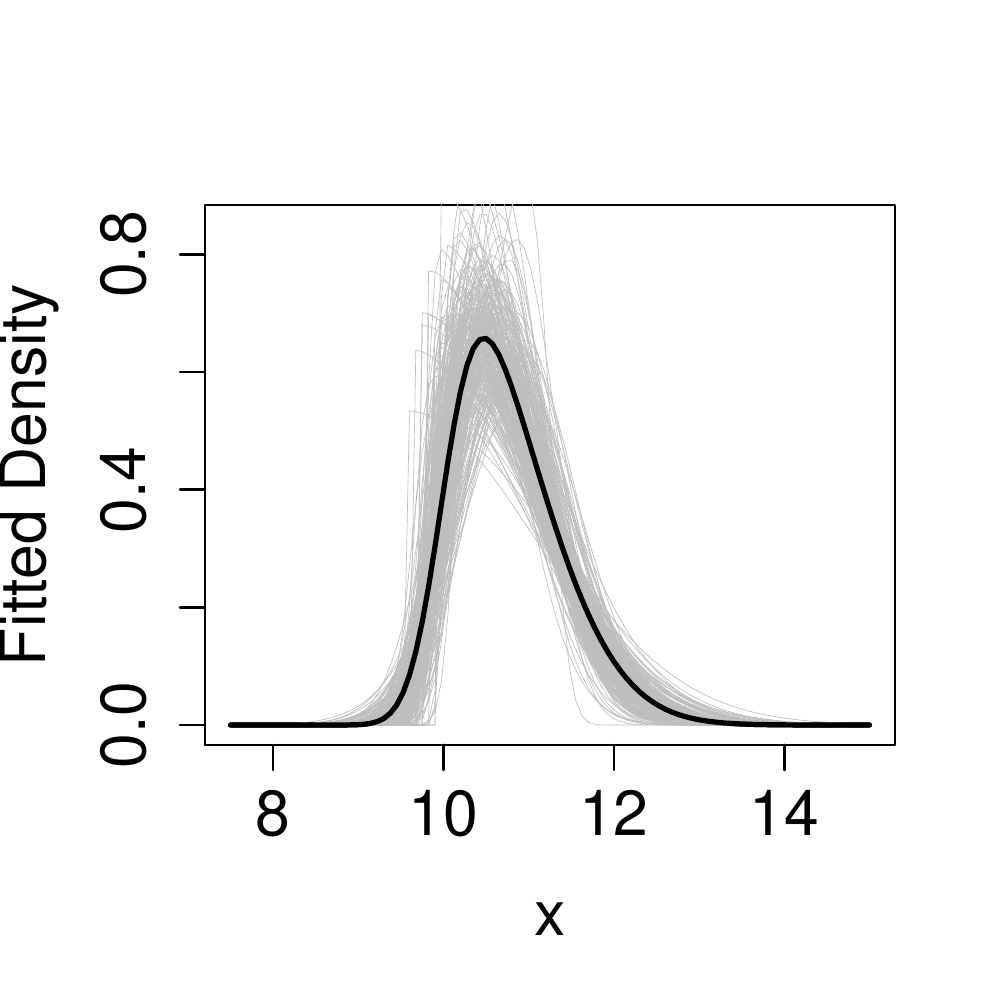} & 
\includegraphics[scale=0.5]{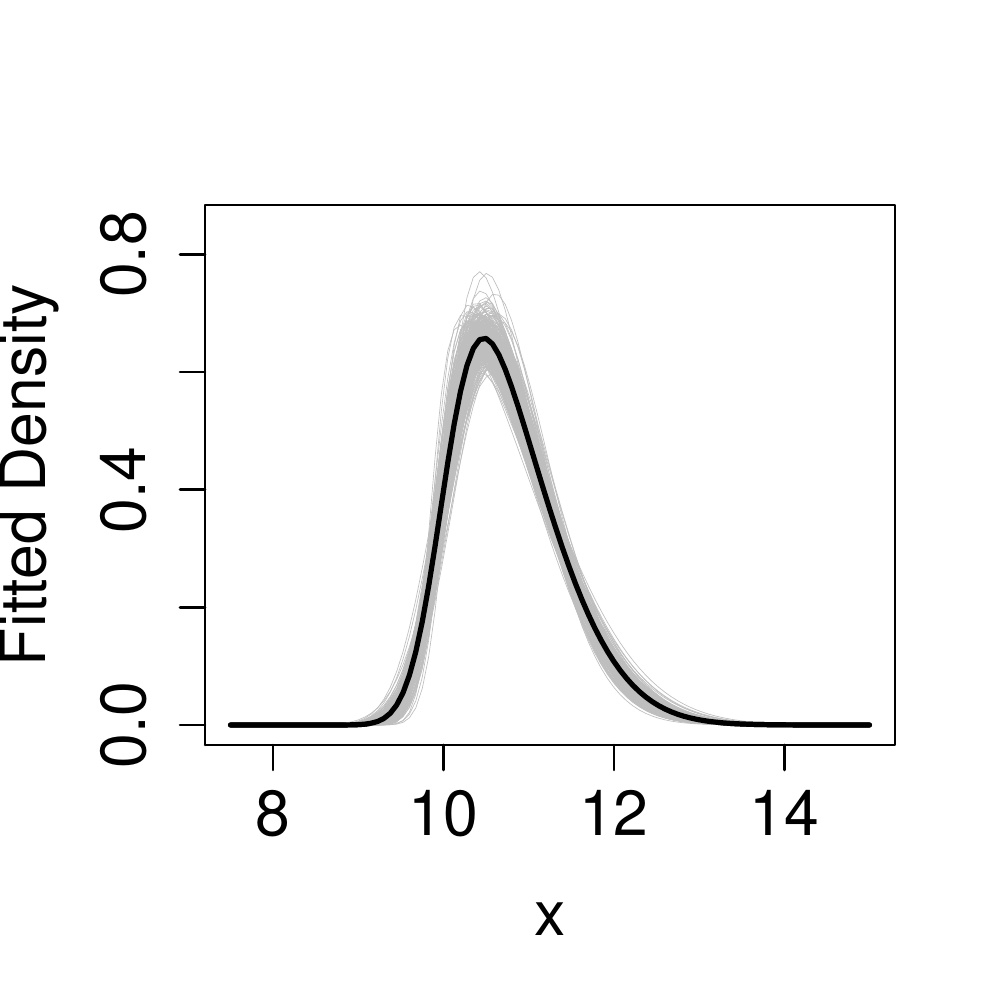} & 
\includegraphics[scale=0.5]{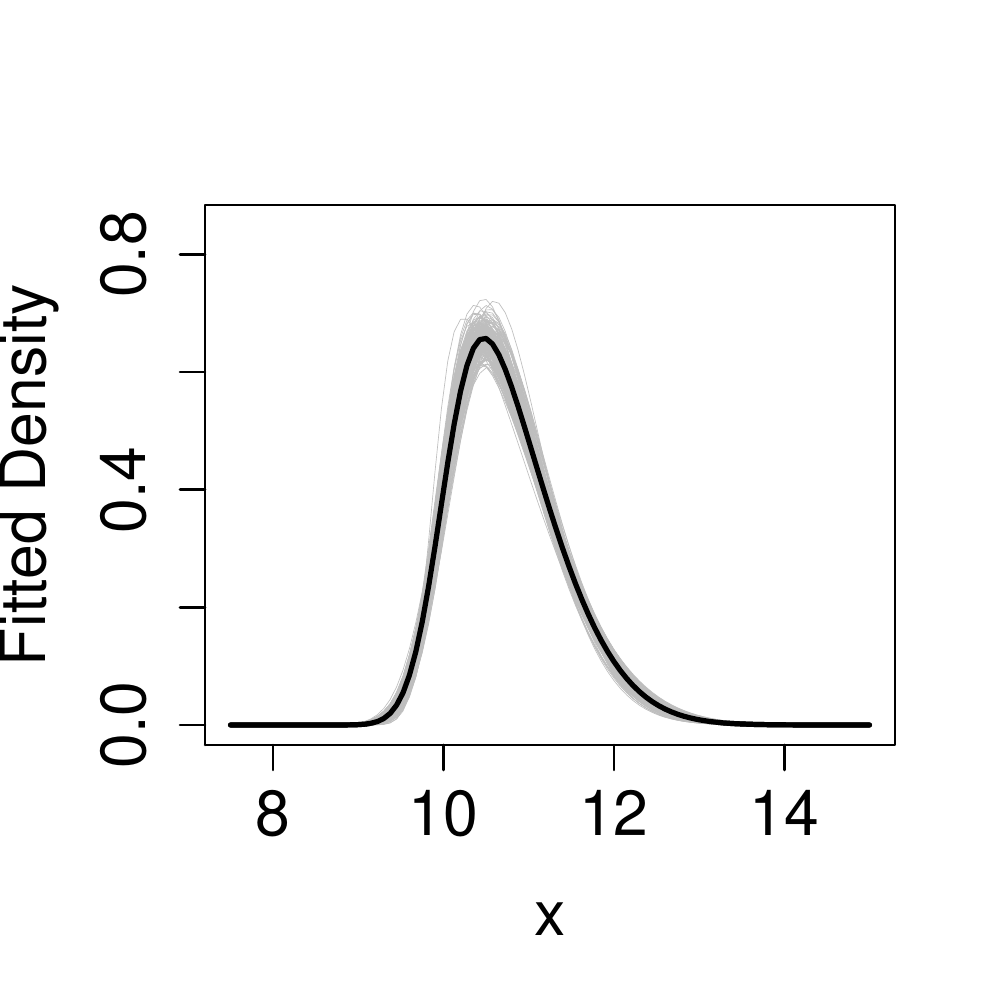} \\
(a) & (b) & (c) \\
\includegraphics[scale=0.5]{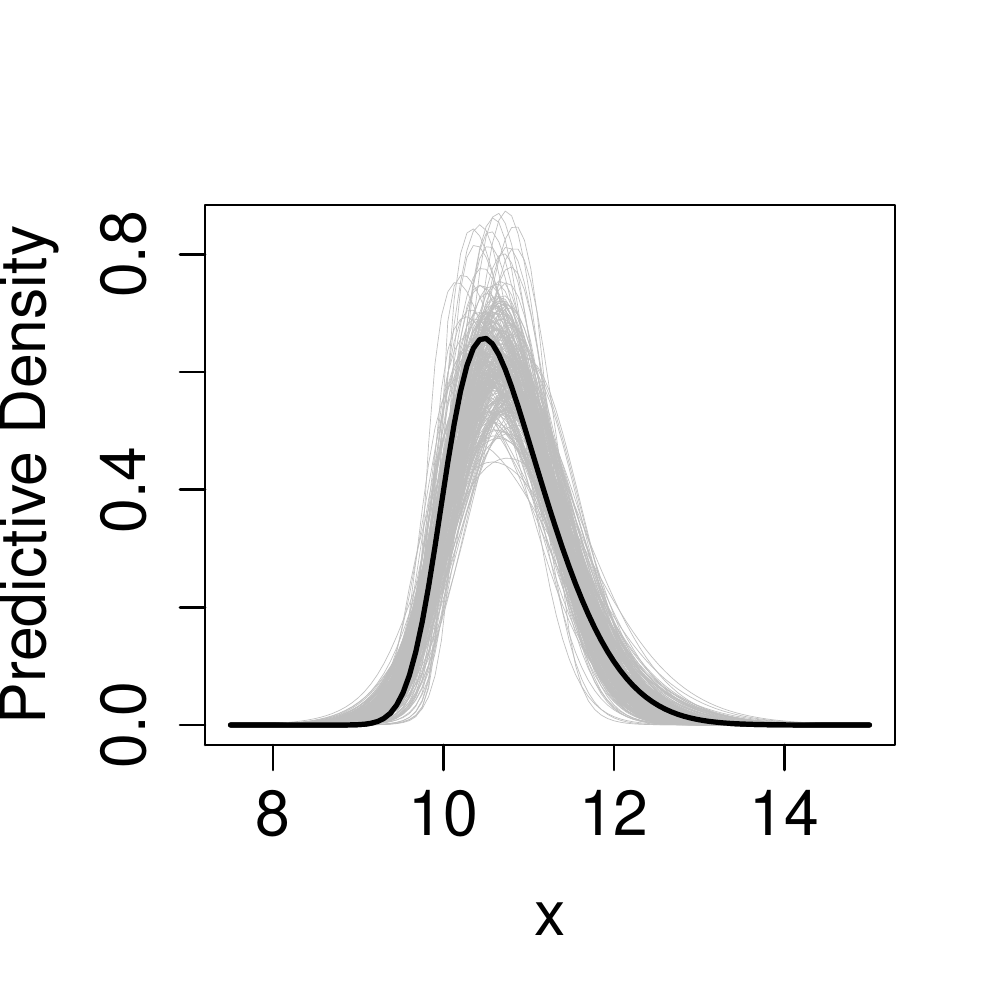} & 
\includegraphics[scale=0.5]{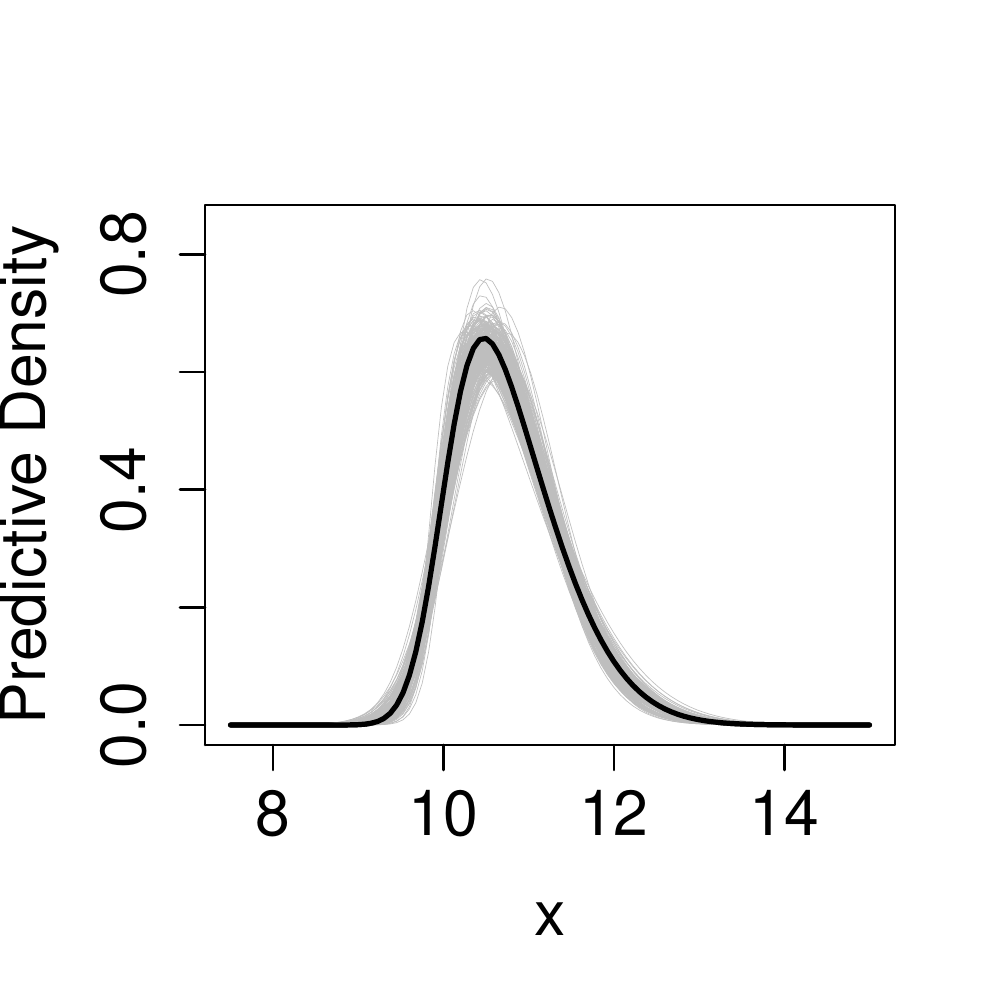} & 
\includegraphics[scale=0.5]{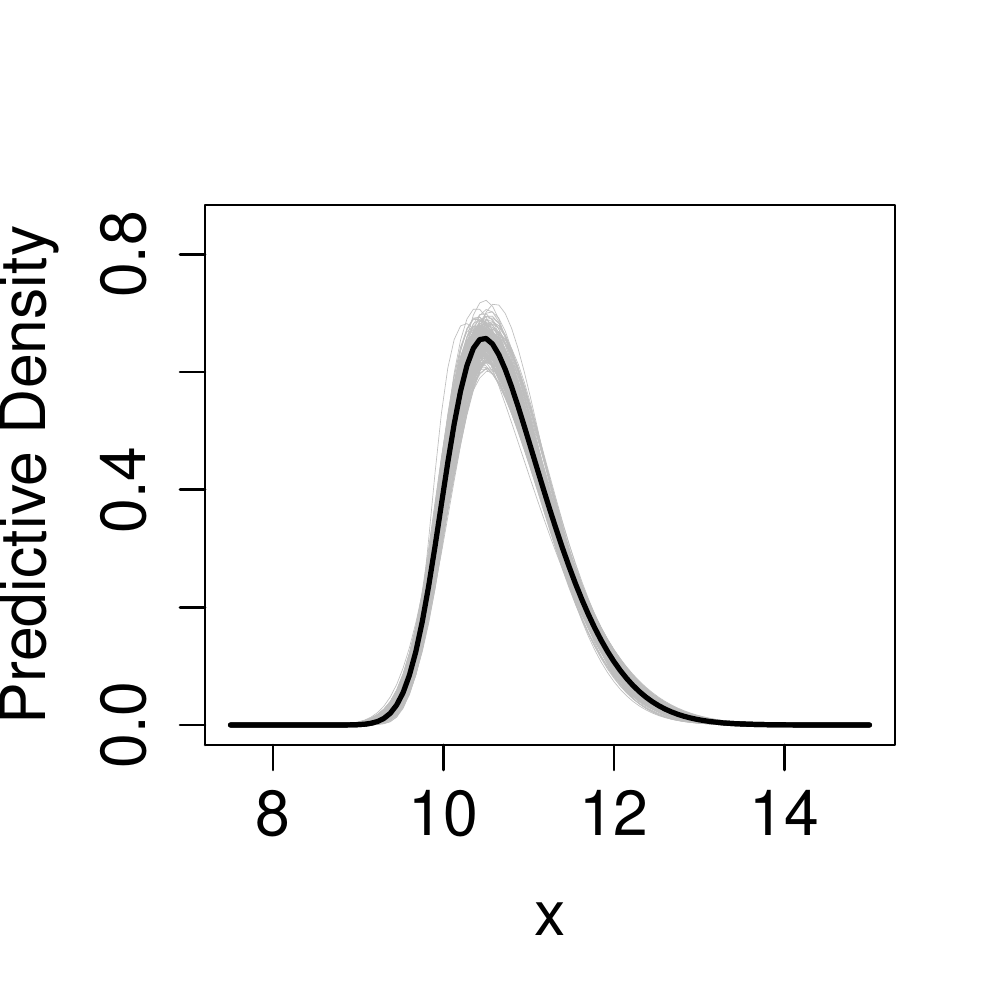} \\
(d) & (e) & (f) \\
\includegraphics[scale=0.5]{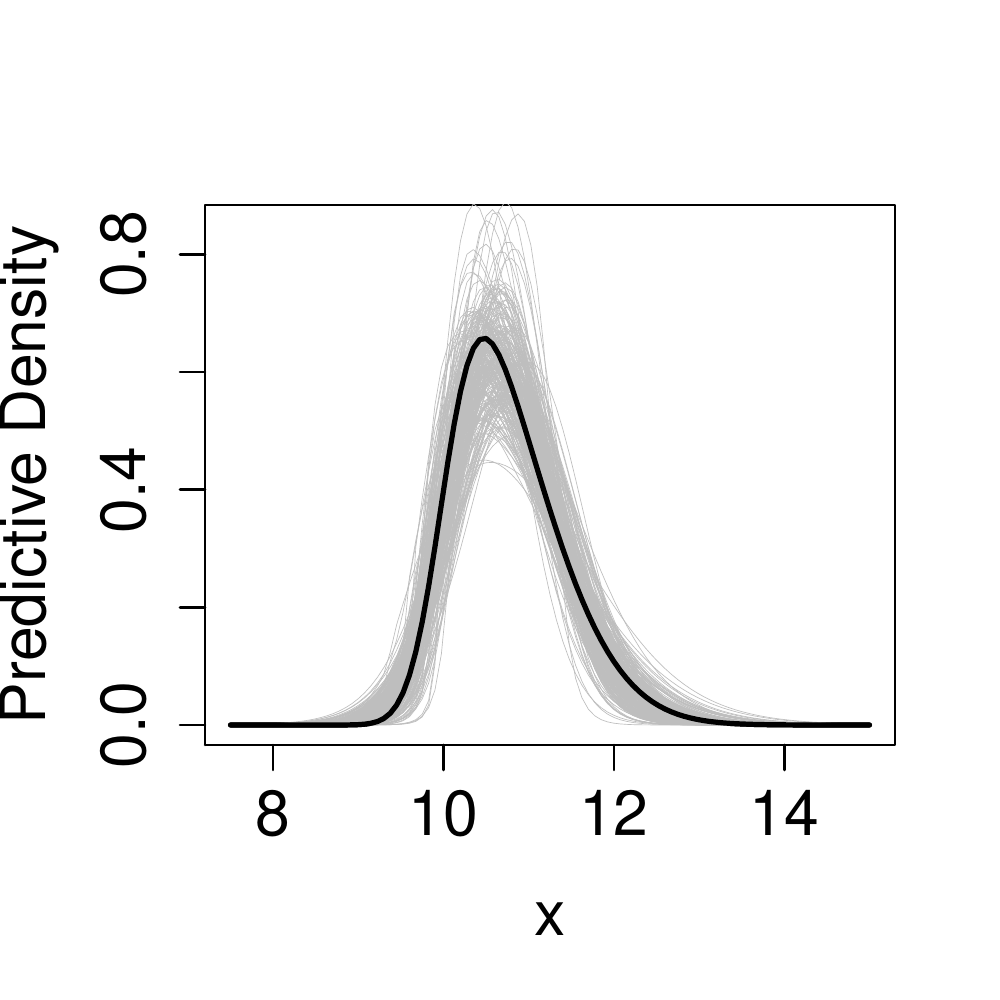} & 
\includegraphics[scale=0.5]{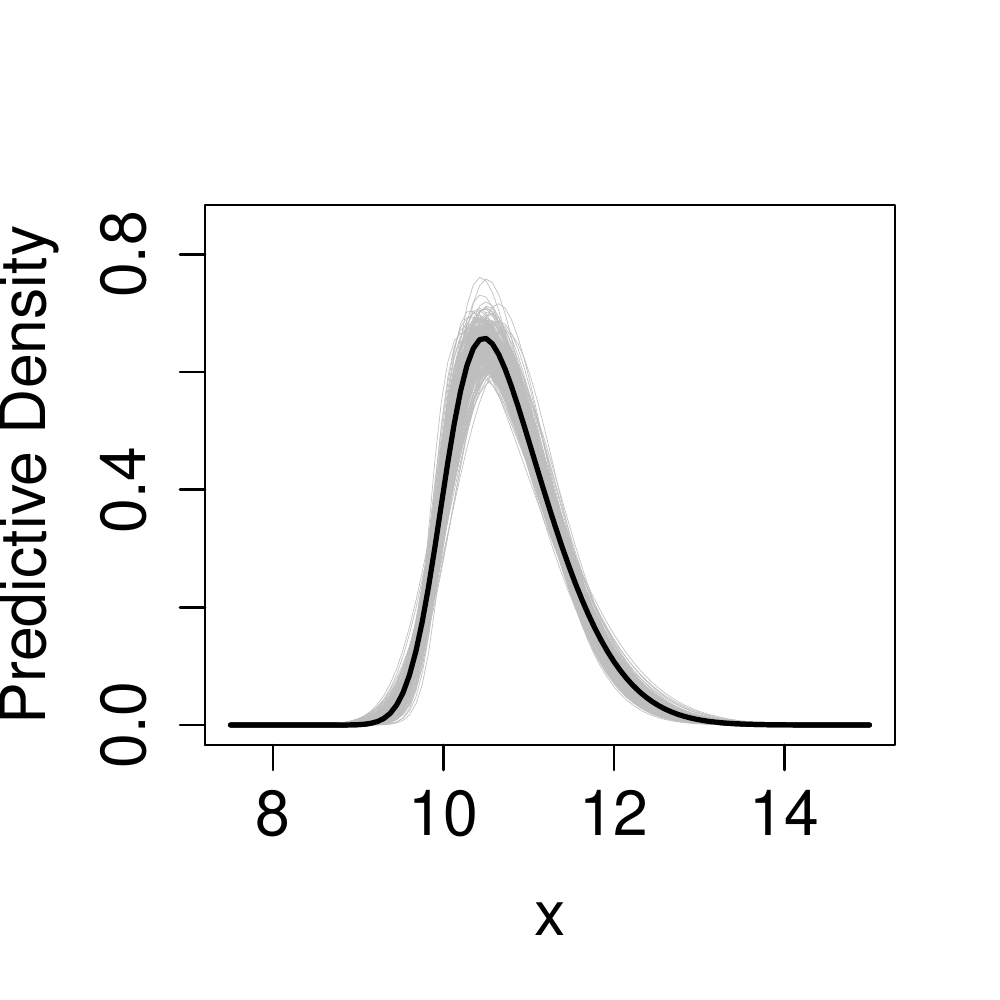} & 
\includegraphics[scale=0.5]{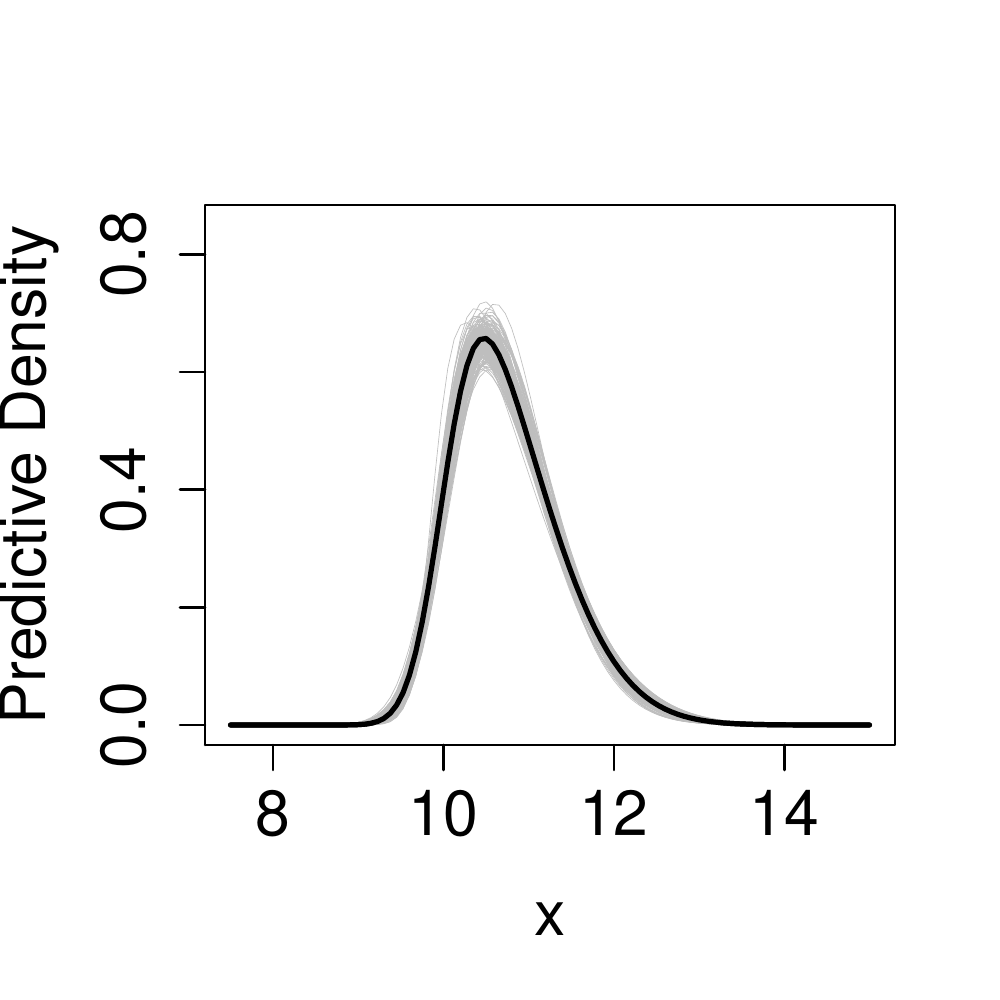} \\
(g) & (h) & (i) \\
\end{tabular}
\caption{(a) MLE-fitted densities ($\lambda = 3$, $n=50$), (b) MLE-fitted densities ($\lambda = 3$, $n=250$), (c) MLE-fitted densities ($\lambda = 3$, $n=500$); (d) Predictive densities for the independence Wasserstein prior ($\lambda = 3$, $n=50$), (e) Predictive densities for the independence Wasserstein prior ($\lambda = 3$, $n=250$), (f) Predictive densities for the independence Wasserstein prior ($\lambda = 3$, $n=500$); (g) Predictive densities for the independence Jeffreys prior ($\lambda = 3$, $n=50$), (d) Predictive densities for the independence Jeffreys prior ($\lambda = 3$, $n=250$), (d) Predictive densities for the independence Jeffreys prior ($\lambda =3$, $n=500$). }
\label{fig:lambda3}
\end{figure}

\begin{figure}
\begin{tabular}{ccc}
\centering
\includegraphics[scale=0.5]{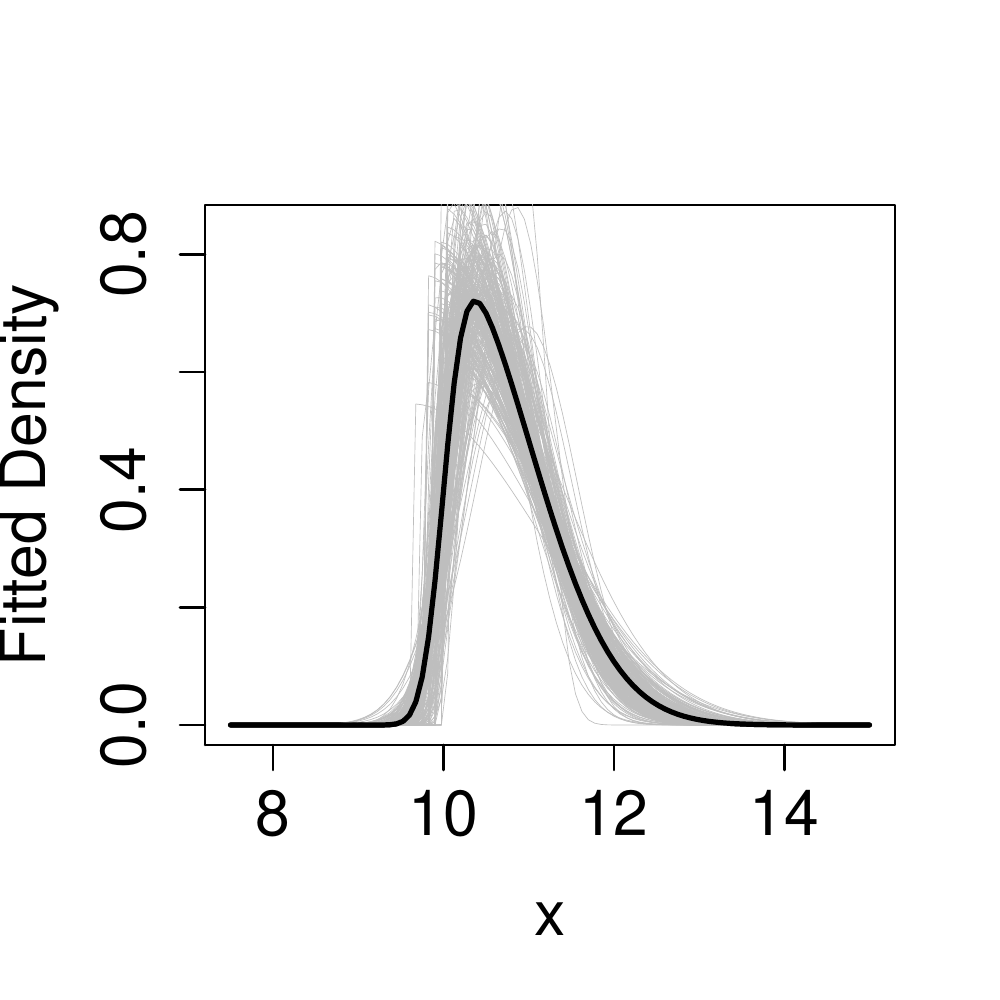} & 
\includegraphics[scale=0.5]{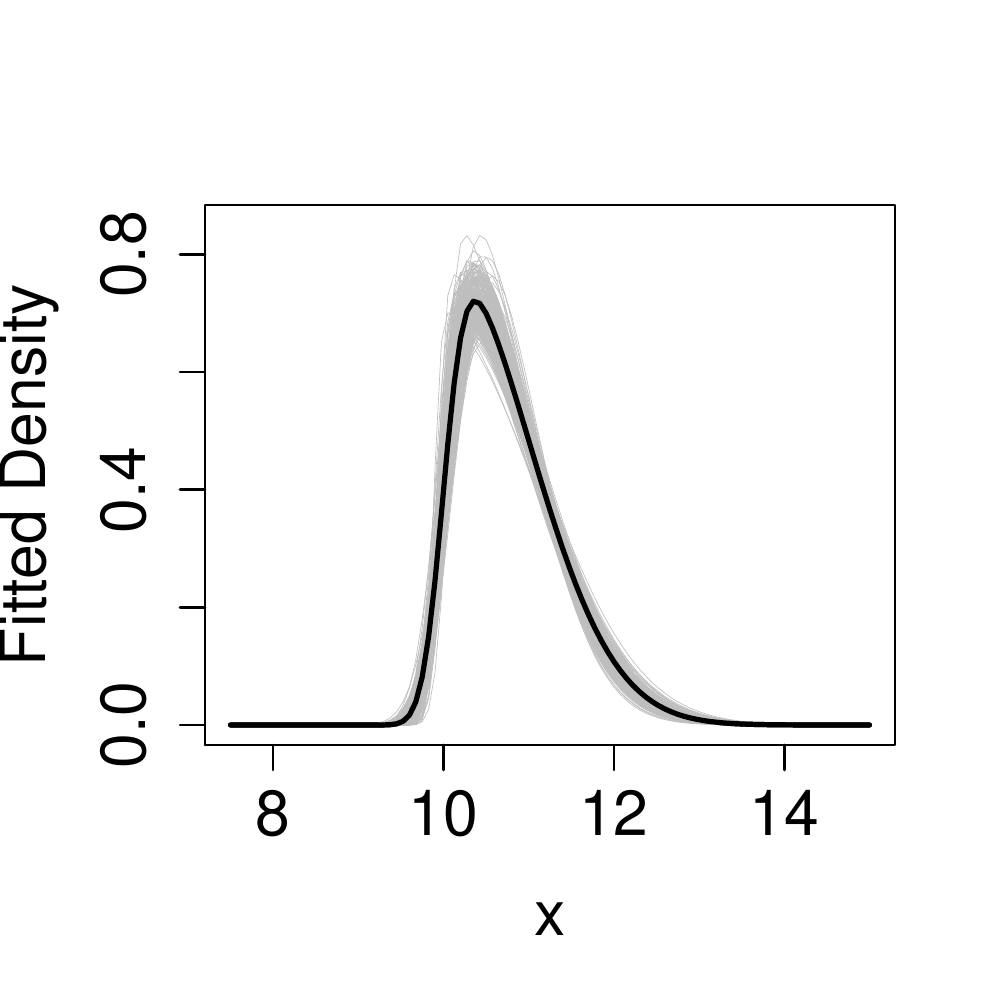} & 
\includegraphics[scale=0.5]{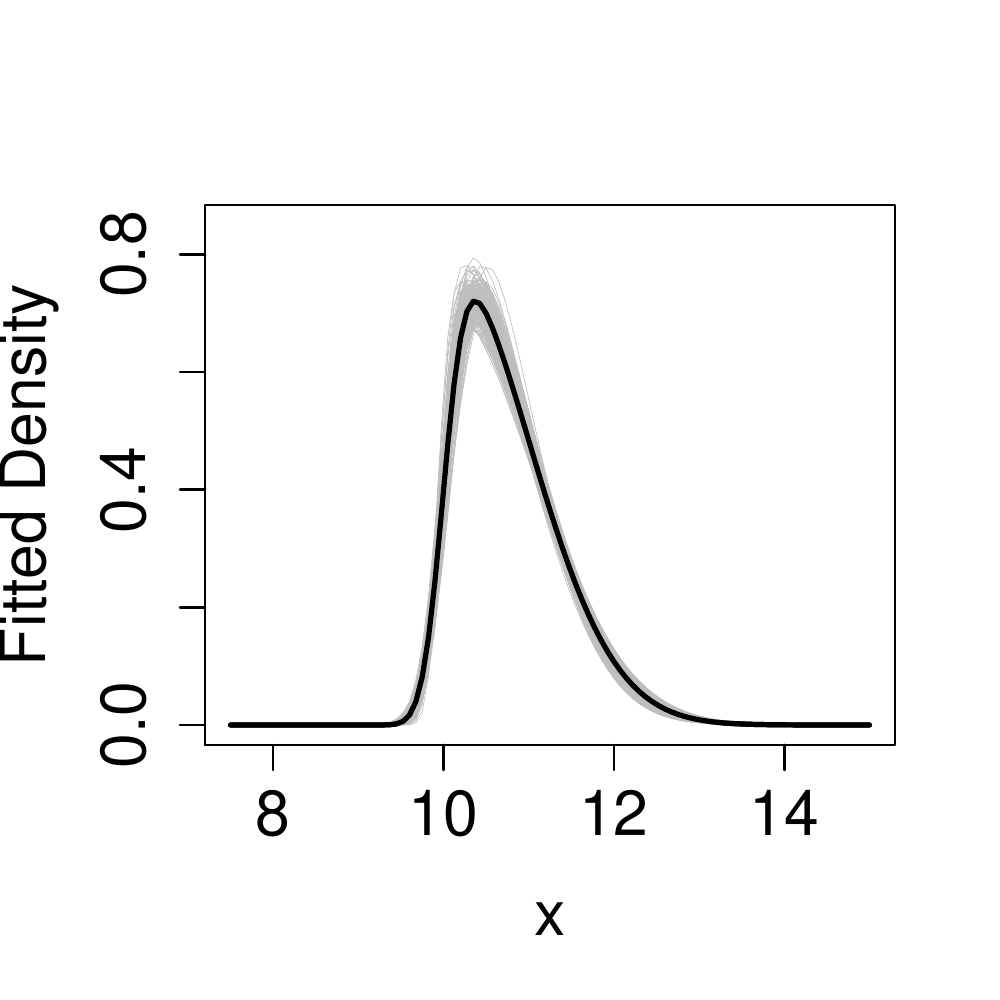} \\
(a) & (b) & (c) \\
\includegraphics[scale=0.5]{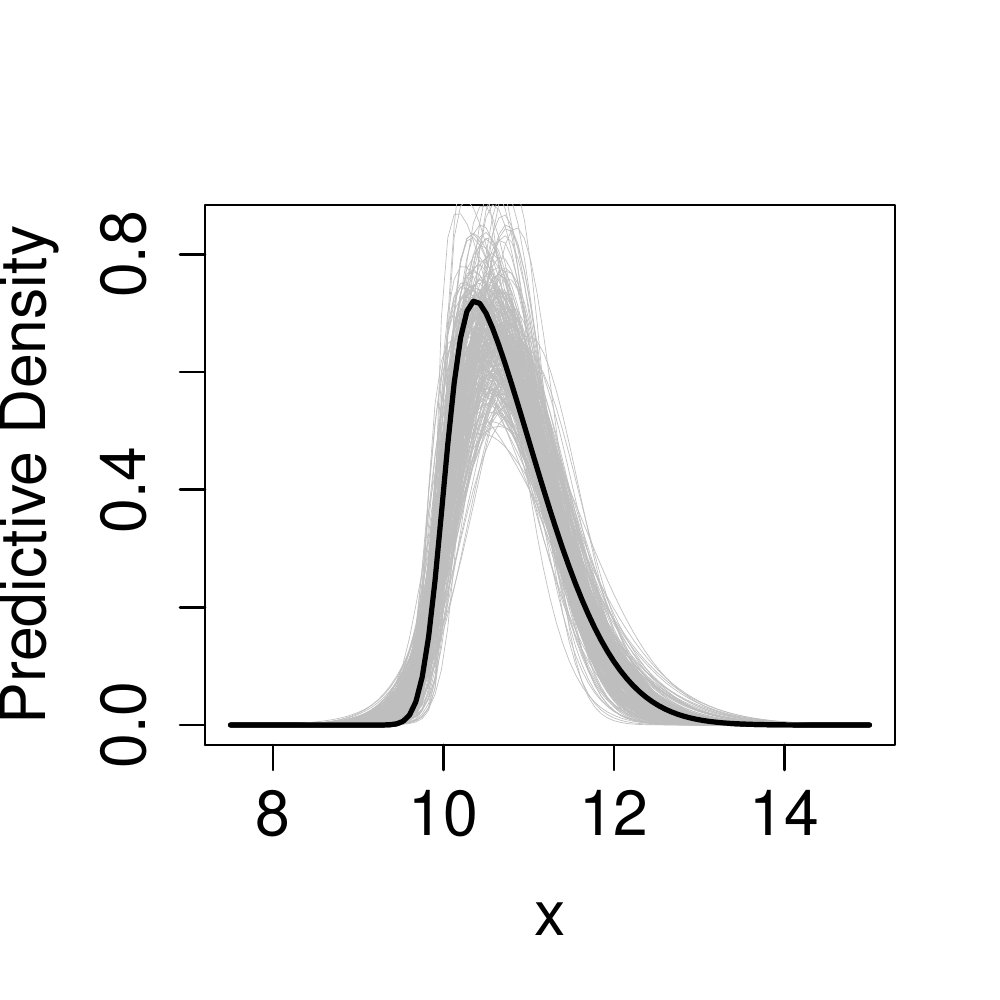} & 
\includegraphics[scale=0.5]{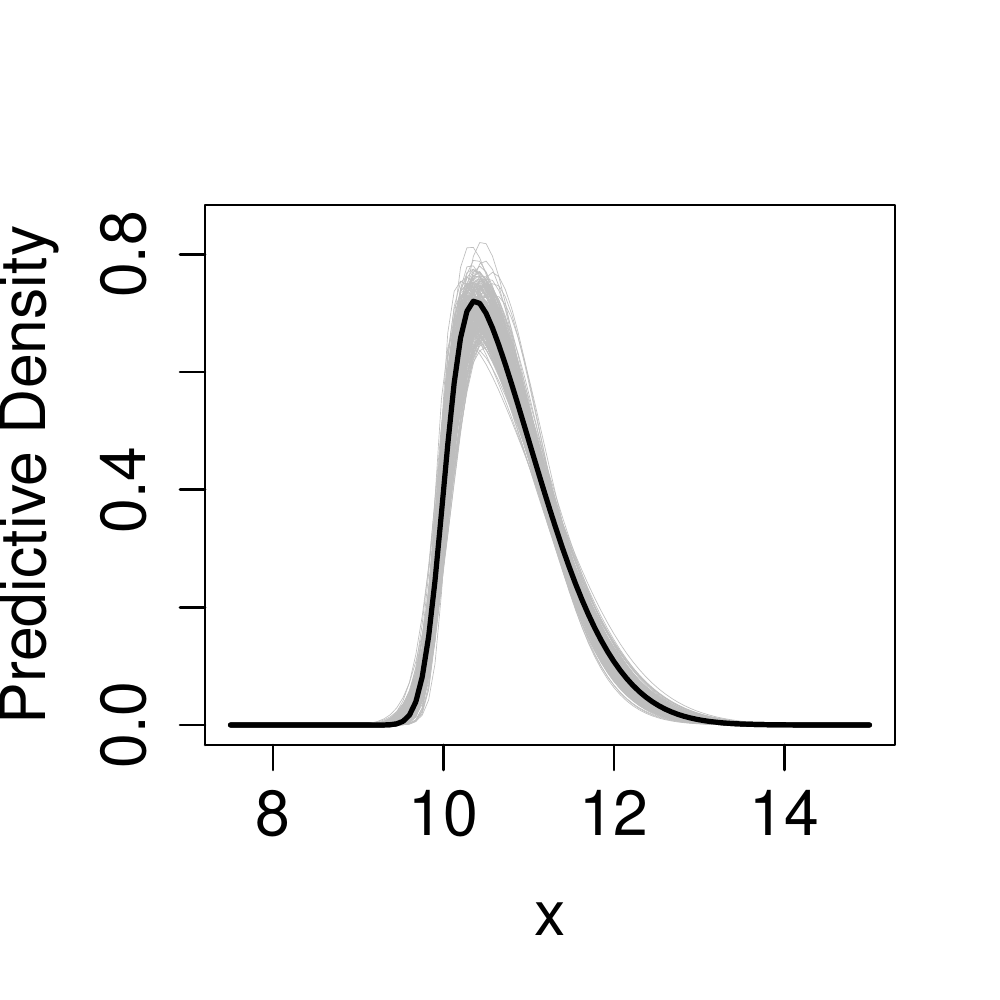} & 
\includegraphics[scale=0.5]{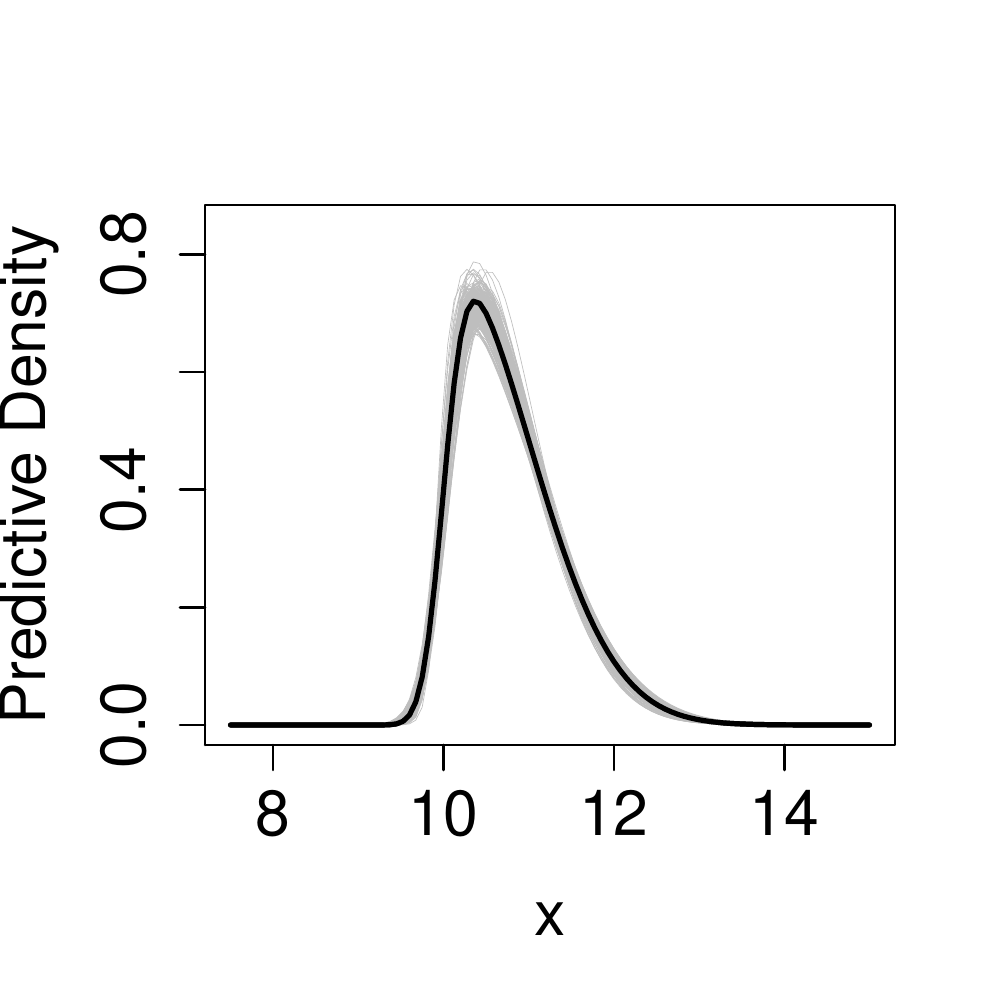} \\
(d) & (e) & (f) \\
\includegraphics[scale=0.5]{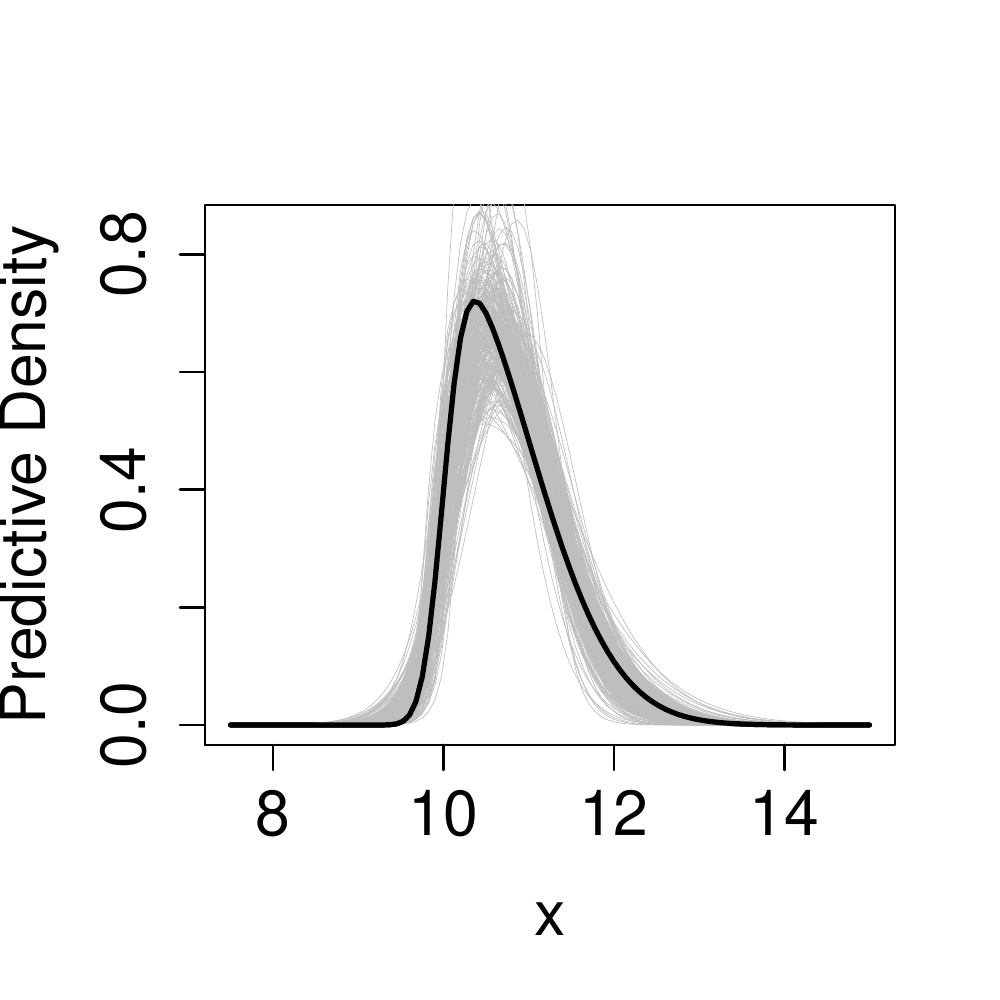} & 
\includegraphics[scale=0.5]{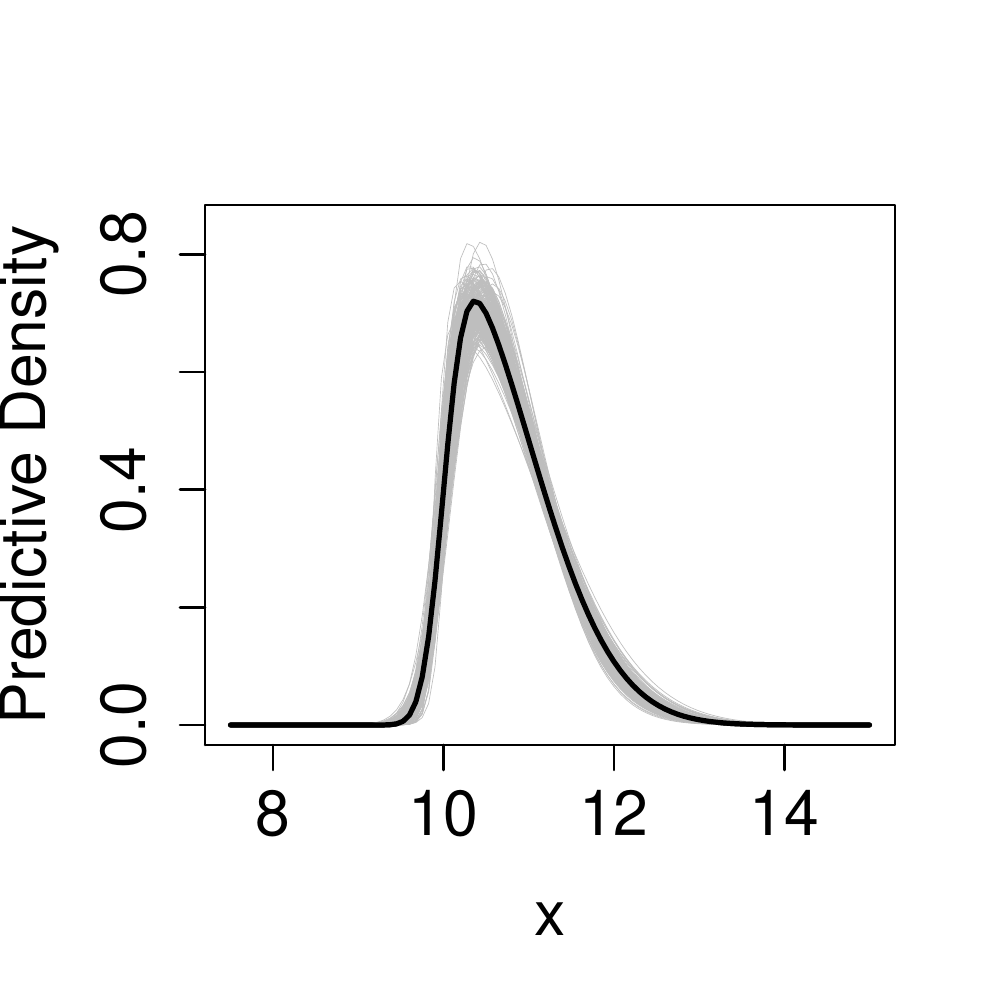} & 
\includegraphics[scale=0.5]{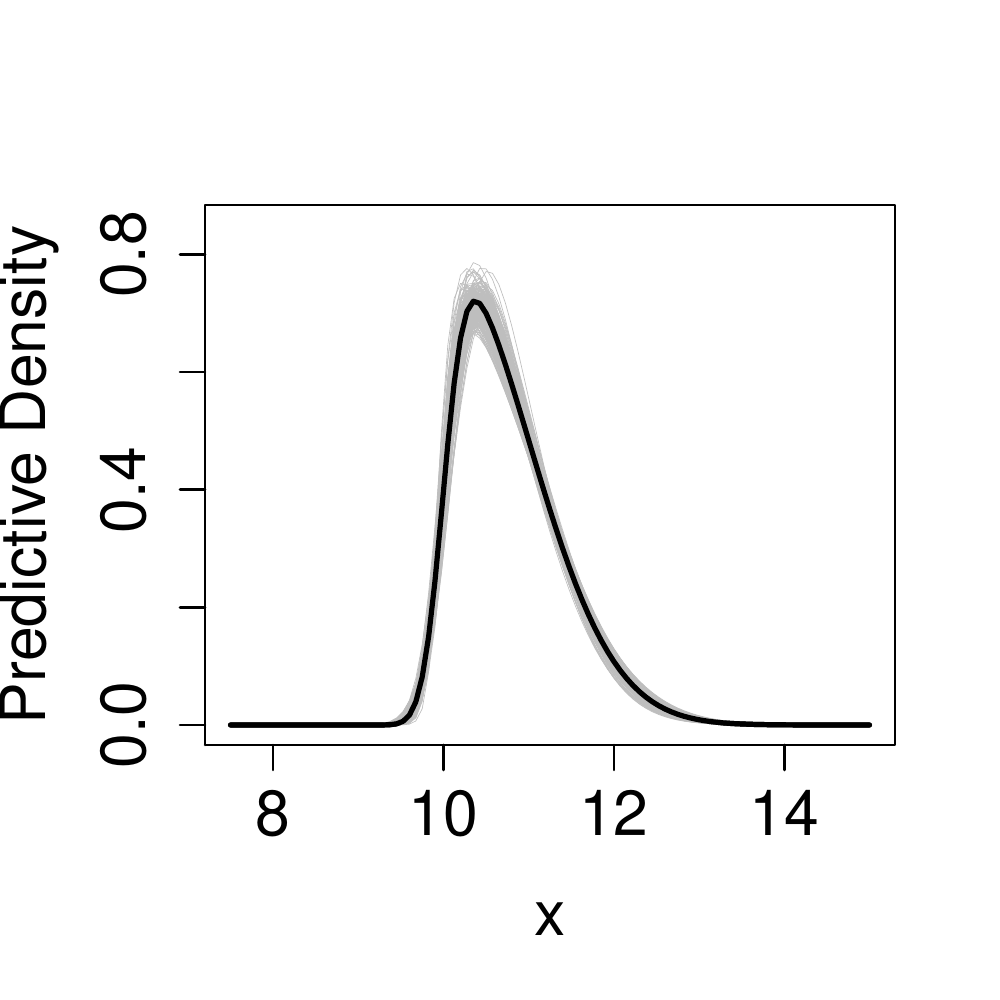} \\
(g) & (h) & (i) \\
\end{tabular}
\caption{(a) MLE-fitted densities ($\lambda = 5$, $n=50$), (b) MLE-fitted densities ($\lambda = 5$, $n=250$), (c) MLE-fitted densities ($\lambda = 5$, $n=500$); (d) Predictive densities for the independence Wasserstein prior ($\lambda = 5$, $n=50$), (e) Predictive densities for the independence Wasserstein prior ($\lambda = 5$, $n=250$), (f) Predictive densities for the independence Wasserstein prior ($\lambda = 5$, $n=500$); (g) Predictive densities for the independence Jeffreys prior ($\lambda = 5$, $n=50$), (d) Predictive densities for the independence Jeffreys prior ($\lambda =5$, $n=250$), (d) Predictive densities for the independence Jeffreys prior ($\lambda =5$, $n=500$). }
\label{fig:lambda5}
\end{figure}

\pagebreak

\begin{table}[ht]
\centering
\begin{tabular}{cccccc}
  \hline
 & $\beta_0$ (1) & $\beta_1$ (0) & $\beta_2$ (0.5) & $\beta_3$ (1) & $\sigma$ (0.5) \\ 
  \hline
  \multicolumn{6}{c}{  $n=50$} \\
mMean & 0.999 & -0.005 & 0.499 & 1.003 & 0.507 \\ 
mSD & 0.077 & 0.134 & 0.097 & 0.122 & 0.055 \\
  mRMSE & 0.102 & 0.184 & 0.127 & 0.162 & 0.073 \\ 
  Coverage & 0.968 & 0.936 & 0.972 & 0.936 & 0.972 \\ 
  mMLE & 0.999 & -0.004 & 0.498 & 1.001 & 0.498 \\ 
  RSME-MLE & 0.074 & 0.138 & 0.089 & 0.116 & 0.050 \\ 
      \multicolumn{6}{c}{  $n=250$} \\
mMean & 0.998 & -0.003 & 0.498 & 1.003 & 0.502 \\ 
  mSD & 0.033 & 0.057 & 0.041 & 0.052 & 0.023 \\
  mRMSE & 0.045 & 0.081 & 0.055 & 0.071 & 0.032 \\ 
  Coverage & 0.948 & 0.924 & 0.956 & 0.932 & 0.932 \\ 
  mMLE & 0.998 & -0.002 & 0.498 & 1.003 & 0.501 \\ 
  RSME-MLE & 0.033 & 0.062 & 0.040 & 0.053 & 0.024 \\
        \multicolumn{6}{c}{  $n=500$} \\
mMean & 1.002 & -0.000 & 0.498 & 0.999 & 0.499 \\ 
  mSD & 0.023 & 0.040 & 0.028 & 0.036 & 0.016 \\ 
  mRMSE & 0.030 & 0.054 & 0.038 & 0.046 & 0.022 \\ 
  Coverage & 0.968 & 0.936 & 0.956 & 0.976 & 0.952 \\ 
  mMLE & 1.002 & 0.000 & 0.498 & 0.998 & 0.498 \\ 
  RSME-MLE & 0.021 & 0.040 & 0.027 & 0.031 & 0.016 \\ 
   \hline
\end{tabular}
\caption{Simulation results for the linear regression model with Wasserstein prior.}
\label{tab:lrm}
\end{table}

\clearpage
\bibliographystyle{plainnat}
\bibliography{references}

\end{document}